\documentclass[11pt]{article}
\usepackage{amssymb}
\usepackage{times}

\title{Uniformisations partielles et crit\`eres \`a la Hurewicz dans le plan.\indent}
\author{Dominique LECOMTE}
\date{\it~Trans. Amer. Math. Soc.\rm ~347, 11 (1995), 4433-4460}

\newcommand{\Ana}{{\it\Sigma}^{1}_{1}}
\newcommand{\Ca}{{\it\Pi}^{1}_{1}}
\newcommand{\Borel}{{\it\Delta}^{1}_{1}}

\newcommand{\boraone}{{\bf\Sigma}^{0}_{1}}
\newcommand{\boratwo}{{\bf\Sigma}^{0}_{2}}
\newcommand{\borathree}{{\bf\Sigma}^{0}_{3}}

\newcommand{\boraxi}{{\bf\Sigma}^{0}_{\xi}}

\newcommand{\borone}{{\bf\Delta}^{0}_{1}}

\newcommand{\bormone}{{\bf\Pi}^{0}_{1}}
\newcommand{\bormtwo}{{\bf\Pi}^{0}_{2}}
\newcommand{\bormthree}{{\bf\Pi}^{0}_{3}}

\newcommand{\bormxi}{{\bf\Pi}^{0}_{\xi}}

\newtheorem{thm} {Th\'eor\`eme} [section]

\newtheorem{defis} [thm] {D\'efinitions}
\newtheorem{cor} [thm] {Corollaire}
\newtheorem{lem} [thm] {Lemme}

\newtheorem{prop} [thm] {Proposition}

\begin{document}

\maketitle

\noindent {\footnotesize {\bf R\'esum\'e.} On donne des caract\'erisations des bor\'eliens potentiellement d'une classe de Wadge donn\'ee, parmi les bor\'eliens \`a coupes verticales d\'enombrables d'un produit de deux espaces polonais. Pour ce faire, on utilise des r\'esultats d'uniformisa-tion partielle.}\bigskip

 Cet article fait suite \`a l'\'etude des classes de Wadge 
potentielles commenc\'ee dans [Le1]. On d\'emontre entre 
autres les r\'esultats annonc\'es dans [Le2]. Je renvoie 
le lecteur \`a [Ku] pour ce qui est des notions de base 
en topologie, et \`a [W] pour ce qui concerne les classes 
de Wadge. On utilisera les notations et notions standard de la 
th\'eorie descriptive des ensembles, ainsi que des notions de th\'eorie 
descriptive effective, qui peuvent \^etre trouv\'ees dans [Mo]. 
Rappelons les d\'efinitions de base :\bigskip

\noindent\bf D\'efinitions.\it\ (a) Soit $\Gamma$ une classe de parties d'espaces 
polonais de dimension 0. On dit que $\Gamma$ est une $classe\ de\ Wadge\ de\ 
bor\acute eliens$ s'il existe un espace polonais $P_0$ de dimension 0, et un bor\'elien $A_0$ de 
$P_0$ tels que pour tout espace polonais $P$ de dimension 0 et pour toute 
partie $A$ de $P$, $A$ est dans $\Gamma$ si et seulement s'il existe une 
fonction continue $f$ de $P$ dans $P_0$ telle que $A = f^{-1}(A_0)$.\smallskip

\noindent (b) Soient $X$ et $Y$ des espaces polonais, et $A$ un bor\'elien de 
$X \times Y$. Si $\Gamma$ est une classe de Wadge, on dira que 
$A$ est $potentiellement\ dans\ \Gamma$ $($ce qu'on notera $A \in 
\mbox{pot}(\Gamma))$ s'il existe des topologies polonaises de dimension 0, 
$\sigma$ $($sur $X)$ et $\tau$ $($sur $Y)$, plus fines que les topologies 
initiales, telles que $A$, consid\'er\'e comme partie de $(X, \sigma) \times 
(Y, \tau)$, soit dans $\Gamma$.\rm\bigskip

 La motivation pour l'\'etude de ces classes de Wadge potentielles trouve son 
origine dans l'\'etude des relations d'\'equivalence bor\'eliennes (ou plus g\'en\'eralement 
des structures bor\'eliennes \`a plusieurs variables). Par exemple,  
dans [HKL], on \'etudie le pr\'e-ordre qui suit. Si $E$ $($resp. $E')$ est 
une relation d'\'equivalence bor\'elienne sur l'espace polonais $X$ $($resp. $X')$, on pose
$$E \leq E' \Leftrightarrow\mbox{ [il existe }f\mbox{ bor\'elienne de }X\mbox{ dans }X'\mbox{ 
telle que }x E y\Leftrightarrow f(x) E' f(y)].$$
En d'autres termes, l'application $f$ d\'efinit par passage au quotient une 
injection de $X/E$ dans $X'/E'$, ce de mani\`ere bor\'elienne. La relation 
$E \leq E'$ peut s'\'ecrire : $E = (f\times f)^{-1} (E')$, avec $f$ bor\'elienne ; 
or si $E'$ est de classe de Wadge $\Gamma$ (ou m\^eme de classe de Wadge 
potentielle $\Gamma$), il est facile de v\'erifier que $E$ est de classe de Wadge 
potentielle $\Gamma$ (alors que $E$ n'est pas de classe $\Gamma$ en g\'en\'eral). La 
classe de Wadge potentielle est donc un invariant naturel du pr\'e-ordre $\leq$. 
Bien que cet invariant soit en g\'en\'eral grossier, il fournit des informations 
sur $\leq$ ; par exemple, dans [Lo1], A. Louveau montre que la classe 
des relations d'\'equivalence $\boraxi$ n'est pas co-finale dans les relations 
d'\'equivalence bor\'eliennes pour le pr\'e-ordre $\leq$, en utilisant la notion de 
classe de Wadge potentielle.

\vfill\eject

 Il en d\'eduit une d\'emonstration simple de la 
non-existence d'une relation d'\'equivalence bor\'elienne maximum pour $\leq$. Ce 
r\'esultat avait par ailleurs \'et\'e d\'emontr\'e ant\'erieurement par H. Friedman et 
L. Stanley, en utilisant les travaux de H. Friedman sur la diagonalisation 
bor\'elienne.\bigskip

  L'un des principaux r\'esultats concernant les classes de Wadge de bor\'eliens 
est l'existence de ``tests d'Hurewicz", dont le principe est : un bor\'elien 
n'est pas d'une classe de Wadge donn\'ee si et seulement s'il est au moins 
aussi compliqu\'e qu'un exemple type n'\'etant pas de cette classe. Hurewicz a 
d\'emontr\'e l'existence du test pour la classe $G_\delta$. Le r\'esultat pr\'ecis est le\bigskip

\noindent\bf Th\'eor\`eme 2.10\it\ Soit $X$ un espace polonais, et $A$ un bor\'elien de $X$. 
Alors $A$ n'est pas $G_\delta$ si et seulement s'il existe $E$ d\'enombrable 
sans point isol\'e tel que $\overline{E}\setminus E \approx \omega^\omega$ et $E = A \cap 
\overline{E}$.\rm\bigskip

 Apr\`es les travaux de nombreux auteurs, parmi lesquels J. R. Steel (voir [S]), 
ainsi que A. Louveau et J. Saint Raymond, l'existence de tests a \'et\'e \'etablie 
pour toutes les classes de Wadge ; dans [Lo-SR], il est d\'emontr\'e :\bigskip

\noindent\bf Th\'eor\`eme.\it\ Si $\xi$ est un ordinal d\'enombrable non nul, il existe un
compact $P_\xi$ de dimension 0 et un vrai $\boraxi$ de $P_\xi$, $A_\xi$, 
tels que si $A$ est un bor\'elien de l'espace polonais $X$, on ait : 
$A$ n'est pas $\bormxi$ de $X$ si et seulement s'il existe $f : P_\xi 
\rightarrow X$ injective continue telle que $A_\xi = f^{-1}(A)$.\rm\bigskip

  L'ensemble $A_\xi$ est dit ``test d'Hurewicz". Un des objectifs de cet 
article est d'\'etudier la possibilit\'e d'obtenir des r\'esultats similaires, pour 
les classes de Wadge potentielles. Dans [Le1], il est d\'emontr\'e un lemme qui sugg\`ere l'int\'er\^et qu'on  
pourrait avoir \`a \'etudier les probl\`emes d'uniformisation partielle, en vue 
de caract\'eriser les ensembles potentiellement ferm\'es. Le test usuel pour 
savoir si un ensemble $A$ est ferm\'e est de prendre une suite convergente de 
points de $A$ et de regarder si la limite est dans $A$. Un tel test ne peut 
pas convenir pour caract\'eriser les ensembles potentiellement ferm\'es, puisqu'un 
singleton peut \^etre rendu ouvert-ferm\'e tout en gardant des topologies polonaises. 
Cependant, on remarque que lorsqu'on raffine la topologie d'un espace polonais, 
tout en gardant une topologie polonaise, les deux topologies co\"\i ncident sur un 
$G_\delta$ dense pour la topologie initiale. D'o\`u l'id\'ee de 
remplacer les points par des graphes de fonctions continues et ouvertes, objets 
qui rencontrent tout produit de deux $G_\delta$ denses si les domaines et images 
des fonctions sont assez ``gros".\bigskip

\noindent\bf Lemme.\it\ Soient $X$ et $Y$ des espaces polonais, $(C_n)~($resp. 
$(D_n))$ des suites d'ouverts non vides de $X$ (resp. $Y$), 
$f_n : C_n \rightarrow D_n$ continues et ouvertes, 
$B := \bigcup\limits_{n\in\omega \setminus \{ 0 \}} \mbox{Gr}(f_n)$, et 
$A$ un bor\'elien de $X \times Y$ contenant $B$; si 
$\overline{B} \setminus A$ contient $\mbox{Gr}(f_0)$, alors $A$ 
est non-$\mbox{pot}(\bormone)$.\rm\bigskip

 La question est de savoir s'il y a une r\'eciproque. Nous allons voir que 
c'est en partie le cas : cette r\'eciproque a lieu, modulo un changement de 
topologies, si $A$ est \`a la fois $\mbox{pot}(\borathree)$ et $\mbox{pot}(\bormthree)$, \`a ceci 
pr\`es qu'on ne peut pas imposer le sens dans lequel se trouvent les graphes. 
Ceci signifie que les fonctions partielles peuvent \^etre d\'efinies sur une 
partie de $X$ et arriver dans $Y$, ou \^etre d\'efinies sur une partie de $Y$ et 
arriver dans $X$. Le r\'esultat pr\'ecis est un cas particulier du th\'eor\`eme 2.3 ; 
si $f_n$ est une fonction partielle de $X$ 
dans $Y$ ou de $Y$ dans $X$, on notera $G_n$ la partie de 
$X\times Y$ \'egale au graphe de $f_n$ si $f_n$ va de $X$ 
dans $Y$, et \'egale \`a $\{(x,y)\in X\times Y~/~x = f_n(y) \}$ sinon.\bigskip
 
\noindent\bf Th\'eor\`eme *.\it\ Soient $X$ et $Y$ des espaces 
polonais, et $A$ un bor\'elien $\mbox{pot}(\borathree)\cap \mbox{pot}(\bormthree)$ de 
$X \times Y$. $A$ est non-$\mbox{pot}(\bormone)$ si et seulement s'il existe 
des espaces polonais $Z$ et $T$ parfaits de dimension 0, des ouverts-ferm\'es 
non vides $A_n$ et $B_n$ (l'un dans $Z$ et l'autre dans $T$), des surjections continues ouvertes $f_n:A_n\rightarrow B_n$, et 
des injections continues $u$ et $v$ tels que $\bigcup_{n\in\omega\setminus\{ 0\}} G_n \subseteq 
(u\times v)^{-1}(A)$, $G_0 \subseteq (u\times v)^{-1}(\check A)$, et 
$G_0 = \overline{\bigcup_{n\in\omega\setminus\{ 0\}} G_n} \setminus (\bigcup_{n\in\omega\setminus\{ 0\}} G_n)$.\rm\bigskip

 Le th\'eor\`eme 2.3 \'etend ce r\'esultat \`a la caract\'erisation des bor\'eliens 
potentiellement diff\'erence transfinie d'ouverts, toujours parmi les bor\'eliens 
$\mbox{pot}(\borathree)\cap \mbox{pot}(\bormthree)$. On ne se ram\`ene donc pas \`a un exemple type, 
comme dans le th\'eor\`eme de A. Louveau et J. Saint Raymond, mais \`a une situation 
type. Pour d\'emontrer le th\'eor\`eme 2.3, l'outil essentiel est le\bigskip

\noindent\bf Th\'eor\`eme 1.13\it\ Soient $X$ et $Y$ des espaces 
polonais parfaits de dimension 0, $A$ un $G_\delta$ l.p.o. 
non vide de $X \times Y$. Alors il existe des ensembles 
presque-ouverts non vides $F$ et $G$, l'un contenu dans $X$ 
et l'autre dans $Y$, et une surjection continue ouverte de 
$F$ sur $G$ dont le graphe est contenu dans $A$ ou dans
$\{(y,x) \in Y \times X~ /~ (x,y) \in A\}$ selon le cas.\rm\bigskip

 Pour le comprendre voici les\bigskip

\noindent\bf D\'efinitions 1.2\it\ (a) Un $G_{\delta}$ d'un espace topologique est 
dit $presque\mbox{-}ouvert$ (ou p.o.) s'il est contenu dans l'int\'erieur de son 
adh\'erence (ce qui revient \`a dire qu'il est dense dans un 
ouvert).\smallskip

\noindent (b) Si $X$ et $Y$ sont des espaces topologiques, une partie $A$
de $X\times Y$ sera dite $localement\ \grave a\ projec\mbox{-}$ $tions\ ouvertes$ (ou l.p.o.) si 
pour tout ouvert $U$ de $X \times Y$, les projections de $A\cap U$ sont 
ouvertes.\rm\bigskip

 Les ensembles l.p.o. se rencontrent par exemple dans la situation suivante : 
$A$ est $\Ana$ dans un produit de deux espaces polonais r\'ecursivement 
pr\'esent\'es. Si on munit ces deux espaces de leur topologie de 
Gandy-Harrington (celle engendr\'ee par les $\Ana$), $A$ devient l.p.o. dans le 
nouveau produit. C'est essentiellement dans cette situation qu'on utilisera 
cette notion, au cours de la section 2.\bigskip

 On \'etudie donc dans un premier temps, et plus largement que n\'ecessaire pour 
la seule \'etude des classes de Wadge potentielles, les probl\`emes 
d'uniformisation partielle, sur des ensembles ``gros" au sens de la cat\'egorie, 
en essayant d'obtenir l'image de la fonction ``grosse" \'egalement. On obtient 
essentiellement des r\'esultats pour les $G_\delta$. Il est \`a noter que malgr\'e 
des hypoth\`eses sym\'etriques, la conclusion du th\'eor\`eme 1.13 ne l'est pas. \bigskip

 Dans [O], il est d\'emontr\'e, sous l'hypoth\`ese du continu, l'existence 
d'une bijection $\Phi$ de $[0,1]$ sur lui m\^eme \'echangeant ensembles maigres et ensembles 
de mesure de Lebesgue nulle. En \'etudiant un exemple o\`u on trouve un graphe 
dans un sens et pas dans l'autre, dans le th\'eor\`eme 1.13, on montre qu'une 
telle application $\Phi$ ne peut pas avoir de propri\'et\'e de mesurabilit\'e 
projective (corollaire 1.8).\bigskip

 Une autre fa\c con de formuler le th\'eor\`eme * est (voir [Le2]) le\bigskip
 
\noindent\bf Th\'eor\`eme.\it\ Soient $X$ et $Y$ des espaces 
polonais, $A$ un bor\'elien de $X \times Y$ \`a coupes verticales d\'enombrables. 
Alors $A$ est non-$\mbox{pot}(\bormone)$ si et seulement s'il 
existe des espaces polonais $Z$ et $T$ parfaits de dimension 0 non vides, des 
fonctions continues $u$ et $v$, des ouverts denses $(A_n)$ de $Z$, 
des applications continues ouvertes $f_n$ de $A_n$ dans $T$, tels que pour tout 
$x$ dans $\bigcap_{n\in\omega} A_n$, $(f_n(x))_{n>0}$ converge simplement 
vers $f_0(x)$, $\mbox{Gr}(g_0) \subseteq 
(u\times v)^{-1}(\check A)$ et 
$\bigcup_{n>0} \mbox{Gr}(g_n) \subseteq (u\times v)^{-1}(A)$.\rm

\vfill\eject

 En d'autres termes, en chaque point du $G_\delta$ dense 
$\bigcap_{n\in\omega} A_n$, on retrouve sur la fibre le test 
usuel pour les ferm\'es. On a un ph\'enom\`ene analogue pour les $G_\delta$ :\bigskip

\noindent\bf Th\'eor\`eme 2.11\it\ Soient $X$ et $Y$ des espaces polonais, 
$A$ un bor\'elien de $X\times Y$ \`a coupes verticales 
d\'enombrables. Alors $A$ est non-$\mbox{pot}(\bormtwo)$ si et 
seulement s'il existe des espaces polonais $Z$ et $T$ parfaits de dimension
 0 non-vides, des injections continues $u$ et 
$v$, des ouverts denses $(A_n)$ de $Z$, des applications 
continues et ouvertes $f_n$ de $A_n$ dans $T$, tels que pour tout $x$ 
dans $\bigcap_{n\in\omega} A_n$, l'ensemble $E_x := \{ 
f_n(x)~/~n\in\omega \}$ soit sans point isol\'e, 
$\overline{E_x} \setminus E_x \approx \omega^\omega$, et 
$E_x = (u\times v)^{-1}(A)_x \cap \overline{E_x}$.\rm

\section{$\!\!\!\!\!\!$ Uniformisation partielle des $G_\delta$.}\indent

 Les probl\`emes d'uniformisation partielle ont d\'ej\`a \'et\'e \'etudi\'es dans [GM], 
o\`u au lieu de consid\'erer la cat\'egorie, il est question
d'ensembles de mesure 1 sur chacun des facteurs. Les r\'esultats qu'on obtient sont \'egalement \`a rapprocher de r\'esultats obtenus par G. Debs et J. Saint Raymond, o\`u il est 
question de fonctions totales et injectives, avec des 
hypoth\`eses de compacit\'e sur chacun des facteurs (cf [D-SR]). Plus pr\'ecis\'ement, il est d\'emontr\'e dans [GM] le r\'esultat suivant :

\begin{thm} Soient $X$ et $Y$ des espaces 
polonais, $\lambda$ (resp. $\mu$) une mesure de 
probabilit\'e sur $X$ (resp. $Y$), et $A$ un bor\'elien de 
$X\times Y$ ayant ses coupes horizontales (resp. 
verticales) non d\'enombrables $\mu$-presque partout (resp. 
$\lambda$-presque partout). Alors il existe un bor\'elien 
$F$ de $X$ (resp. $G$ de $Y$) tels que $\lambda (F) = 
\mu (G) = 1$, et un isomorphisme bor\'elien de $F$ sur $G$ dont 
le graphe est contenu dans $A$.\end{thm}

 On peut se demander si on a un r\'esultat analogue en 
rempla\c cant ``ensemble de mesure 0" par ``ensemble maigre". 
On va voir que non. Pour ce faire on montre un lemme que 
nous r\'eutiliserons.

\begin{defis} (a) Un $G_\delta$ d'un espace topologique est dit $presque\mbox{-}ouvert$ (ou p.o.) s'il est contenu dans l'int\'erieur de sn adh\'erence (ce qui revient \`a dire qu'il est dense dans un ouvert).\smallskip

\noindent (b) Si $X$ et $Y$ sont des espaces topologiques, une partie $A$ de $X\times Y$ sera dite $localement\ \grave a\ projec\mbox{-}$ $tions\ ouvertes$ (ou l.p.o.) si pour tout ouvert $U$ de $X\times Y$, les projections de $A\cap U$ sont ouvertes.\end{defis}
 
\begin{lem} Il existe un espace $X$ polonais parfait 
de dimension 0, et un ferm\'e $A$ l.p.o. \`a coupes parfaites 
non vides de $X^2$, tels que si $F$ et $G$ sont 
presque-ouverts non vides dans $X$ et $f:F\rightarrow G$ 
surjective continue ouverte, $\mbox{Gr}(f) \not \subseteq A$.\end{lem}

\noindent\bf D\'emonstration.\rm\ Soient $X = \omega^\omega$, $\phi$ un hom\'eomorphisme 
de $X$ sur l'ensemble $P_{\infty}$ des suites de 0 et de 1 ayant une infinit\'e de 1, et $\psi$ un hom\'eomorphisme 
de $X$ sur l'ensemble ${\cal K}_P(2^\omega)$ des compacts 
parfaits non vides de $2^\omega$ ; $\psi$ existe car cet ensemble 
ne contient pas les compacts finis, est dense et est 
$G_{\delta}$ de ${\cal K}(2^\omega)\setminus \{ \emptyset \}$ : 
en effet, si on d\'esigne par $N(n,2^\omega)$ le 
$n^{\mbox{i\`eme}}$ ouvert-ferm\'e de base de $2^\omega$, on a 
$$K\ \mbox{est parfait}\Leftrightarrow\!\left\{\!\!\!\!\!\!\! 
\begin{array}{ll}
& \forall m \in \omega~[(K\cap N(m,2^\omega) = \emptyset)\ \mbox{ou}~(\exists~(p,q) \in 
\omega^2~~N(p,2^\omega)\cap N(q,2^\omega) = \emptyset\ \mbox{et}\cr & \cr
& N(p,2^\omega)\cup N(q,2^\omega)\!\subseteq\! N(m,2^\omega)\ \mbox{et}~
N(p,2^\omega)\cap K \not= \emptyset\ \mbox{et}~N(q,2^\omega)\cap K 
\not= \emptyset )].
\end{array}\right.$$

\vfill\eject

 Posons $A := (\phi \times \psi)^{-1}(A')$, o\`u 
$A' := \{ (x,K) \in P_{\infty}\times {\cal K}_P(2^\omega)~/~x\in K \}$. 
Il suffit de montrer les propri\'et\'es pour $A'$. Il est ferm\'e :
$$x\notin K \Leftrightarrow \exists~n\in\omega~~\exists~U\in 
\boraone\lceil2^\omega~ ~N(n,2^\omega) \cap U = \emptyset\ ~\mbox{et}~~
x\in  N(n,2^\omega)\ ~\mbox{et}~~K \subseteq U.$$ 
Il est clairement \`a coupes parfaites non vides, et il est l.p.o., car si 
$N_s:=\{\alpha\in 2^\omega/s\prec\alpha\}$ et 
$$W_{U,V_0,...,V_n}:=\{ K\in {\cal K}_P(2^\omega)/K\subseteq U\ \mbox{et}~\forall i\leq n\ V_i\cap K\not=\emptyset\}$$ 
est l'ouvert-ferm\'e de base de ${\cal K}_P(2^\omega)$, on a 
$$\Pi_1[A'\cap (N_s \times W_{U,V_0,...,V_n})] =\emptyset ~~\mbox{ou}~~N_s \cap U \cap P_\infty
\mbox{, ~et}$$
$$\Pi_2[A'\cap (N_s \times W_{U,V_0,...,V_n})] = \emptyset\ ~\mbox{ou}~\ W_{U,V_0,...,V_n,N_s}\cap 
{\cal K}_P(2^\omega).$$
Raisonnons par l'absurde : il existe $f: F \rightarrow G$ 
surjective continue ouverte dont le graphe est contenu 
dans $A'$, et des ouverts non vides $\cal U$ et $\cal V$ tels que 
$F$ (resp. $G$) soit $G_\delta$ dense de $\cal U$ (resp. $\cal V$).\bigskip

 L'ouvert $\cal V$ \'etant non vide contient un ouvert-ferm\'e 
de base $W_{U,V_0,...,V_n}$, avec $U$, $V_0,~ ...,~V_n$ 
ouverts-ferm\'es tels que les $U \cap V_i$ soient non vides. 
De sorte qu'on peut supposer, quitte \`a se restreindre \`a 
son image r\'eciproque, que ${\cal V} = W_{U,V_0,...,V_n}$.\bigskip

 Soit $(S_j)_{j<n+3}$ une partition de $U$ en 
ouverts-ferm\'es telle que $\forall~j < n + 3$, $\forall~i\leq n$, 
$S_j\cap V_i \not= \emptyset$. Posons, si $I \in {\cal P}(n + 3) \setminus \{\emptyset \}$,
$$O_I := \{K \in W_{U,V_0,...,V_n}~ /~ I = \{j < n + 3~/~K \cap S_j \not= \emptyset\}\}.$$
Alors $(O_I)$ est une partition en ouverts-ferm\'es non 
vides de $W_{U,V_0,...,V_n}$, donc $(f^{-1}(O_I))$ est 
une partition en ouverts-ferm\'es non vides de $F$ ; par 
propri\'et\'e de r\'eduction, on trouve une suite d'ouverts 
non vides de $\cal U$ deux \`a deux disjoints, $(W_I)$, telle que 
$f^{-1}(O_I) = F\cap W_I$.\bigskip

 Alors si $j < n+3$, $W_{\{j\}} \subseteq S_j$, sinon on 
trouve $x$ dans $F\cap W_{\{j\}} \setminus S_j$, et $x\in 
f(x)\in  O_{\{j\}}$, donc $f(x) \subseteq S_j$, une 
contradiction.\bigskip

 Maintenant si $I \subseteq n + 3$ est de cardinal 2, il 
existe $i \leq n$ tel que $V_i \cap \bigcup_{j\in I} S_j \subseteq 
\overline{W_I}$. En effet, par ce qui pr\'ec\`ede, si $j \in I$, 
$S_j \not\subseteq \overline{W_I}$ sinon $W_{\{j\}} \cap W_I\not= \emptyset$. 
Donc si pour tout $i \leq n$ il existe $j$ dans $I$ tel que 
$V_i \cap S_j \not\subseteq \overline{W_I}$, l'ensemble 
$\{K \in O_I ~/~ K \subseteq \check{\overline{W_I}}\}$ est 
ouvert non vide de $\overline{G}$, donc rencontre $G$ en $K = f(x)$ 
et $x \in F \cap W_I$, donc $K \cap W_I$ est non vide, une 
contradiction.\bigskip

 En particulier, si $1 \leq j \leq n + 2$, il existe $i_j \leq n$ 
tel que $V_{i_j} \cap (S_0 \cup S_j) \subseteq 
\overline{W_{\{0,j\}}}$. Soient $j \not= j'$ tels que $i_j = 
i_{j'}$ ; alors $V_{i_j} \cap S_0 \subseteq 
\overline{W_{\{0,j\}}} 
\cap \overline{W_{\{0,j'\}}}$, ce qui contredit la disjonction 
de $W_{\{0,j\}}$ et $W_{\{0,j'\}}$.$\hfill\square$

\begin{thm} Soient $X := {\cal K}(2^\omega) \setminus 
\{\emptyset \}$, $Y := 2^\omega$, et 
$A := \{(K,x) \in X \times Y~ /~ x\in  K\}$ ; alors $A$ 
est un ferm\'e \`a coupes horizontales non d\'enombrables, \`a 
coupes verticales non d\'enombrables sur un ensemble 
co-maigre de $X$, mais si $F$ (resp. $G$) est un bor\'elien 
co-maigre de $X$ (resp. $Y$), il n'existe aucun 
isomorphisme bor\'elien de $F$ sur $G$ dont le graphe soit 
contenu dans $A$.\end{thm}

\noindent\bf D\'emonstration.\rm\ Comme on l'a vu dans la preuve du lemme 
pr\'ec\'edent, $A$ est ferm\'e dans $X \times Y$, et l'ensemble 
des compacts parfaits non vides est $G_\delta$ dense de $X$.\bigskip

 L'ensemble des compacts non d\'enombrables est donc 
co-maigre dans $X$, et les coupes verticales de $A$ sont 
donc non d\'enombrables, sauf sur un maigre.\bigskip

 Raisonnons par l'absurde : il existe des ensembles $F$ et 
$G$, ainsi qu'un isomorphisme bor\'elien $f$ comme dans 
l'\'enonc\'e.\bigskip

 Alors il existe un $G_\delta$ dense $F'$ (resp. $G'$) de 
$X$ (resp. $Y$) tels que $F' \subseteq {\cal K}_P(2^\omega)\cap F$, 
$G' \subseteq G \cap P_\infty$, et que 
$f: F' \rightarrow G'$ soit un hom\'eomorphisme.\bigskip

 En effet, on remarque que si $U$ est un ouvert dense de $Y$, 
alors $\{K \in X ~/~K\subseteq U\}$ est un ouvert dense de 
$X$ ; par suite, si $M$ est maigre relativement \`a $G$, on 
a l'inclusion 
$$f^{-1}(M) \subseteq \Pi_X[(X \times M) \cap A] = 
\{K \in X ~/~ K \cap M \not=\emptyset\}.$$ 
Donc $f^{-1}(M)$ est maigre relativement \`a $X$, et aussi relativement \`a $F$ 
qui est co-maigre dans $X$.\bigskip

 Soit $G_1$ un $G_\delta$ dense de $Y$, contenu dans 
$G \cap P_\infty$, sur lequel $f^{-1}$ est continue ; ce 
qui pr\'ec\`ede montre que $f^{-1}(G_1)$ est co-maigre dans $X$.\bigskip

 Si $(U_n)$ est une base de la topologie de $G_1$, 
$f^{-1}(U_n)$ est bor\'elien de $X$, donc \'egal \`a un ouvert 
$V_n$ modulo un maigre $M_n$ ; choisissons pour $F'$ un 
$G_\delta$ dense de $X$ contenu dans 
$${\cal K}_P(2^\omega)\cap F \cap f^{-1}(G_1) \setminus (\bigcup_{n\in\omega} M_n).$$
Il reste \`a poser $G' := f''F'$.\bigskip

 En effet, on a $f^{-1}(U_n \cap G') = F' \cap f^{-1}(U_n) 
= F'\cap V_n$ ; $G' = (f^{-1})^{-1}(F') \cap G_1$ est 
$G_\delta$. Et si $G'$ n'\'etait pas dense dans $Y$, on 
trouverait un ouvert non vide $U$ de $Y$ disjoint de $G'$ ; 
mais $\{K \in X ~/~ K \subseteq U\}$ serait un ouvert non 
vide, et rencontrerait donc $F'$ en un point $K$ qui 
v\'erifierait $f(K) \in K \subseteq U$ et aussi $f(K) \in G'$, 
ce qui est la contradiction cherch\'ee.\bigskip

 Mais ceci contredit la preuve du lemme 1.3.$\hfill\square$\bigskip

 En analysant les raisons de ce r\'esultat n\'egatif, nous 
allons maintenant montrer que l'application transformant 
ensembles de mesure 0 en ensembles maigres n'a pas de 
propri\'et\'es de mesurabilit\'e projective, sous hypoth\`ese de 
d\'etermination des jeux projectifs.\bigskip

 Pour d\'emontrer le th\'eor\`eme 1.1, les auteurs montrent le 
lemme suivant :

\begin{lem} Soient $X$ un espace polonais, $\mu$ 
une mesure de probabilit\'e sur $X$, et $A$ un bor\'elien de 
$X \times X$ ayant ses coupes horizontales non 
d\'enombrables $\mu$-presque partout. Alors il existe un 
$K_\sigma$ de $X$, $F$, et une application bor\'elienne 
$f: F \rightarrow X$ tels que $\mbox{Gr}(f) \subseteq A$ et 
$$\mu(\{y \in X ~/~ f^{-1}(\{y\})\ \mbox{est non d\'enombrable}\}) = 1.$$\end{lem}

\begin{cor} Soient $B := [0,1]$, $\lambda$ 
la mesure de Lebesgue sur $B$, $A$ un bor\'elien de 
$B \times B$ tel que $\lambda(\{y \in B ~/~ A^y$ est non 
d\'enombrable$\}) = 1$ ; alors il existe des bor\'eliens 
disjoints de $B$, $B_0$ et $B_1$, tels que $\lambda(\{y \in B ~/~
B_i \cap A^y$ est non d\'enombrable$\}) = 1$.\end{cor}

\noindent\bf D\'emonstration.\rm\ Soient $F$ et $f$ fournis par le 
lemme 1.5, et $G_y$ une copie de $2^\omega$ contenue dans 
$f^{-1}(\{y\})$, si ce bor\'elien est non d\'enombrable. Si 
$x\in G_y \setminus \{\mbox{min}~G_y,~\mbox{max}~ G_y\}$, les ouverts 
$[0,x[\cap G_y$ et $]x,1] \cap G_y$ sont non vides, donc 
non d\'enombrables ; posons donc 
$$E(y,x)\ \Leftrightarrow\ x\in F\ \mbox{et}~f(x) = y\ \mbox{et}~
[0,x[\cap f^{-1}(\{y\})\mbox{,}~]x ,1] \cap f^{-1}(\{y\})\ \mbox{sont non d\'enombrables.}$$
Alors $E$ est analytique dans $B \times B$, donc par le 
th\'eor\`eme de von Neumann on trouve une application 
$g: B \rightarrow B$ Baire-mesurable telle que si 
$f^{-1}(\{y\})$ est non d\'enombrable, $g(y) \in F$, 
$f(g(y)) = y$, et $[0,g(y)[\cap f^{-1}(\{y\})$, 
$]g(y),1] \cap f^{-1}(\{y\})$ sont non d\'enombrables.\bigskip

 Soit alors $G$ un bor\'elien de $B$ tel que $\lambda(G) = 1$, 
la restriction de $g$ \`a $G$ soit bor\'elienne, et contenu 
dans $\{y\in  B ~/~ f^{-1}(\{y\})$ est non d\'enombrable$\}$. 
Il reste \`a poser\bigskip

 \centerline{$B_0 := \{x \in F ~/~ f(x) \in G$ et $x < g(f(x))\}$, 
$B_1 := \{x \in F ~/~ f(x) \in G$ et $x > g(f(x))\}.$}\bigskip

 En effet, si $y\in  G$, $[0,g(y)[\cap f^{-1}(\{y\}) \subseteq 
B_0 \cap A^y$, donc $B_0 \cap A^y$ est non d\'enombrable, de 
m\^eme que $B_1\cap A^y$.$\hfill\square$

\begin{lem} Il existe un $G_\delta$ de 
$B \times B$, $A$, o\`u $B := [0,1]$, tel que $\{y \in B ~/~ A^y\ \mbox{est non d\'enombrable}\}$ 
soit co-maigre, alors que pour tout couple $(B_0,  B_1)$ de bor\'eliens disjoints de $B$, 
$$\{y \in B ~/~ B_i \cap A^y\ \mbox{est non d\'enombrable}\}$$ 
ne sont pas tous deux co-maigres.\end{lem}

\noindent\bf D\'emonstration.\rm\ $[0,1] \setminus \mathbb{Q}$~ est hom\'eomorphe \`a 
$P_\infty$, qui est co-d\'enombrable dans l'espace parfait 
$2^\omega$. Il suffit donc de trouver $A$ dans $2^\omega\times2^\omega$, 
ayant les propri\'et\'es du lemme. D'autre part, ${\cal K}(2^\omega) \setminus \{\emptyset\}$ est un 
espace compact m\'etrisable parfait de dimension 0 et non 
vide, donc est hom\'eomorphe \`a $2^\omega$ ; il suffit donc de 
trouver $A$ dans $2^\omega\times {\cal K}(2^\omega)\setminus 
\{\emptyset \}$. Reprenons l'exemple de 1.4 : 
$$A = \{(x,K)\in2^\omega\times {\cal K}(2^\omega)\setminus 
\{\emptyset \}~/~x \in K \}.$$ 
Raisonnons par l'absurde : $B_0$ et $B_1$ existent. Par la preuve de 1.4, ces deux ensembles sont 
n\'ecessairement non maigres ; on trouve donc des ouverts 
disjoints et non vides de $2^\omega$, $U_0$ et $U_1$, tels 
que $B_0\Delta U_0$ et $\check B_0\Delta U_1$ soient 
maigres.\bigskip

 Comme on l'a vu dans la preuve du th\'eor\`eme 1.4, 
$\{K \in {\cal K}(2^\omega)\setminus \{\emptyset 
\}~ /~ K\cap (B_0 \Delta U_0) = \emptyset\}$ est 
co-maigre, ainsi que $\{K \in {\cal K}(2^\omega)\setminus \{\emptyset 
\}~/~B_0\cap K$ est non d\'enombrable$\}$, donc que 
l'ensemble $\{K \in {\cal K}(2^\omega)\setminus \{\emptyset 
\}~/~ K \cap U_0 \not= \emptyset\}$ ; ce dernier 
rencontre donc l'ouvert non vide 
$$\{K \in {\cal K}(2^\omega)\setminus \{\emptyset \}~/~K \subseteq U_1 \}\mbox{,}$$ 
ce qui contredit la disjonction de $U_0$ et $U_1$.$\hfill\square$\bigskip

 Il r\'esulte imm\'ediatement du th\'eor\`eme 19.6 de [O] 
que sous l'hypoth\`ese du continu, il existe une bijection 
idempotente $\Phi$ de $B := [0, 1]$ sur lui-m\^eme telle que $E$ 
est maigre si et seulement si $\lambda (\Phi^{-1}(E)) = 0$. 
Ce qui pr\'ec\`ede entra\^\i ne le 

\begin{cor} Une telle fonction $\Phi$ n'est 
pas bor\'elienne. De plus, sous hypoth\`ese de d\'etermination 
des jeux ${\bf \Delta}^1_{2n+3}$,~\rm $\Phi$ n'est pas 
${\bf \Pi}^1_{2n+1}\rm$-mesurable.\end{cor}

\noindent\bf D\'emonstration.\rm\ Le premier point r\'esulte aussit\^ot du 
corollaire 1.6 et du lemme 1.7. Si $\Phi$ \'etait 
${\bf \Pi}^1_{2n+1}\rm$-mesurable, et si $A$ est le $G_\delta$ 
du lemme 1.7, l'ensemble $A' := (I \times \Phi)^{-1}(A)$ 
serait ${\bf \Pi}^1_{2n+1}\rm$, ainsi que 
$E := \{(y, K) \in B \times {\cal K}(B) \setminus \{\emptyset\} 
~/~K$ est parfait et $K \subseteq A'^y\}$ ; par la 
d\'etermination des jeux ${\bf \Delta}^1_{2n}\rm$, $E$ serait 
uniformisable par un graphe partiel ${\bf \Pi}^1_{2n+1}\rm$,  
et la fonction correspondante se prolongerait en une 
fonction $f$ totale ${\bf \Delta}^1_{2n+3}\rm$-mesurable. Par 
la d\'etermination des jeux ${\bf \Delta}^1_{2n+3}\rm$, $f$ 
serait $\lambda$-mesurable, donc il existerait un bor\'elien 
$G$ de $B$ tel que $\lambda(G) = 1$, la restriction de $f$ 
\`a $G$ soit bor\'elienne, et $(y,f(y)) \in E$ si $y\in G$.\bigskip

 Posons $A'' := \{(x,y) \in B \times B ~/~ y \in G$ et 
$x\in  f(y)\}$ ; alors $A''$ serait un bor\'elien contenu 
dans $A'$, et le corollaire 1.6 fournirait des bor\'eliens 
disjoints $B_0$ et $B_1$. Les ensembles 
$$\{y \in B ~/~ B_i\cap (I \times\Phi)[A'']^y\ \mbox{est non d\'enombrable}\}$$ 
seraient donc co-maigres, ainsi que $\{y \in B ~/~ B_i \cap 
A^y$ est non d\'enombrable$\}$, ce qui contredirait le lemme 
1.7.$\hfill\square$\bigskip

 Dans [Ma], la question suivante est pos\'ee : \'etant 
donn\'e un bor\'elien $A$ de $[0,1] \times [0,1]$ dont toutes 
les coupes sont non maigres, existe-t-il un isomorphisme 
bor\'elien de $[0,1]$ sur lui-m\^eme dont le graphe soit contenu dans $A$ ? La r\'eponse 
est non, comme le montre le th\'eor\`eme 3 de [D-SR]. 
On a cependant le r\'esultat suivant :

\begin{thm} Soient $X$ et $Y$ des espaces 
polonais parfaits non vides, $A \subseteq X \times Y$ ayant 
la propri\'et\'e de Baire et ses coupes non maigres (sauf sur 
des ensembles maigres). Alors il existe un ensemble $K_\sigma \cup 
G_\delta$, $F$ (resp. $G$), co-maigre dans $X$ (resp. $Y$), 
et un isomorphisme de deuxi\`eme classe de $F$ sur $G$ dont 
le graphe est contenu dans $A$.\end{thm}

\begin{lem} Soient $(C_n)$, $(D_n)$ des suites 
d'ouverts-ferm\'es non vides de $\omega^\omega$, $A$ dans $\omega^\omega\times \omega^\omega$ 
tel que $A \Delta (\bigcup_{n\in\omega} C_n \times D_n)$ 
soit maigre ; alors il existe des suites $(K^0_n)$ et $(K^1_n)$ 
de copies de $2^\omega$, des suites $(G^0_n)$ et $(G^1_n)$ 
de $G_\delta$ de $\omega^\omega$ v\'erifiant\smallskip

\noindent (i) $(K^0_n \times G^0_n) \cup (G^1_n \times K^1_n)\subseteq A$\smallskip

\noindent (ii) $K^0_n \subseteq C_n \setminus (\bigcup_{q<n} K^0_q),~
K^1_n \subseteq D_n \setminus (\bigcup_{q<n} K^1_q)$\smallskip

\noindent (iii) $G^0_n$ (resp. $G^1_n$) est dense dans $D_n \setminus 
(\bigcup_{q<n} D_q \cup \bigcup_{p\in\omega} K^1_p)~($resp. 
$C_n \setminus (\bigcup_{q<n} C_q \cup 
\bigcup_{p\in\omega} K^0_p))$\end{lem}

\noindent\bf D\'emonstration.\rm\ On construit d'abord la suite $(K^0_n)$ 
et une suite $(H^0_n)$ de $G_\delta$ de $\omega^\omega$ v\'erifiant 
$$K^0_n \times H^0_n \subseteq A \cap [C_n \setminus 
(\bigcup_{q<n} K^0_q) \times D_n \setminus (\bigcup_{q<n} 
D_q)]\mbox{,}$$
$H^0_n$ \'etant dense dans $D'_n := D_n \setminus 
(\bigcup_{q<n} D_q)$.\bigskip

 Admettons la construction effectu\'ee pour $q < n$. Soit 
$C'_n := C_n\setminus (\bigcup_{q<n} K^0_q)$ ; 
$A\cap C'_n\times D'_n$ est co-maigre dans $C'_n \times D'_n$, 
d'o\`u l'existence de $K^0_n$ et $H^0_n$.\bigskip

 De m\^eme, on construit des suites $(K^1_n)$ et $(H^1_n)$ 
v\'erifiant des propri\'et\'es analogues : 
$$H^1_n \times K^1_n \subseteq A \cap [C_n \setminus 
(\bigcup_{q<n} C_q) \times D_n \setminus (\bigcup_{q<n} 
K^1_q)]\mbox{,}$$ 
et $H^1_n$ est dense dans $C_n \setminus 
(\bigcup_{q<n} C_q)$. Il reste \`a poser 
$G^0_n := H^0_n \setminus 
(\bigcup_{p\in\omega} K^1_p)$, $G^1_n := H^1_n \setminus 
(\bigcup_{p\in\omega} K^0_p).$
Ceci termine la preuve.$\hfill\square$\bigskip 

\noindent\bf D\'emonstration du th\'eor\`eme 1.9.\rm\ L'ensemble $A$ ayant la 
propri\'et\'e de Baire, il est r\'eunion disjointe d'un $G_\delta$ 
et d'un ensemble maigre, qui a ses coupes maigres, sauf 
sur des ensembles maigres ; on peut donc supposer que $A$ 
est $G_\delta$ \`a coupes non maigres de $\omega^\omega\times\omega^\omega$, 
$X$ et $Y$ \'etant parfaits non vides.\bigskip

 On trouve alors des suites $(C_n)$ et $(D_n)$ v\'erifiant 
les conditions du lemme 1.10, qui peut donc s'appliquer. 
Posons $I := \{n\in\omega~/~G^0_n \not= \emptyset \}$, 
$F_0 := \bigcup_{n\in I} K^0_n$, 
$G_0 := \bigcup_{n\in I} G^0_n$. Si $n \in\omega$, $D_n 
\setminus (\bigcup_{q<n} D_q \cup \bigcup_{p\in\omega} K^1_p)$, 
donc $G^0_n$, sont denses dans 
$D_n \setminus (\bigcup_{q<n} D_q)$ ; par cons\'equent, $G_0$ est dense dans $\bigcup_{n\in\omega} D_n$. D'autre 
part, $G_0 \in \borone\mbox{-PU}(G_\delta)\lceil(\bigcup_{n\in\omega} D_n)$, 
donc $G_0$ est $G_\delta$. Comme $A$ est \`a coupes non 
maigres, $\bigcup_{n\in\omega} D_n$ est dense dans 
$Y$, donc $G_0$ est $G_\delta$ dense dans $Y$.\bigskip

 Si $n \in I$, $G^0_n$ est polonais non d\'enombrable, donc 
il existe un isomorphisme de premi\`ere classe $f_n$ de $K^0_n$ 
sur $G^0_n$. Alors 
$$g_0 : \left\{\!\!
\begin{array}{ll} F_0 
& \!\!\!\rightarrow G_0 \cr x 
& \!\!\!\mapsto f_n(x)~\ {\rm si}~\ x\in K^0_n
\end{array}\right.$$ 
est bien d\'efinie (par (ii)), injective (par (iii)) donc bijective, de premi\`ere classe par compacit\'e 
de $K^0_n$, et $g_0^{-1}$ est aussi de premi\`ere classe car 
$G^0_n = G_0 \cap [D_n \setminus (\bigcup_{q<n} D_q)]$. 
Enfin, le graphe de $g_0$ est contenu dans $A$ par (i).\bigskip

 De m\^eme, on trouve un isomorphisme de premi\`ere classe $g_1$, 
du $K_\sigma$ maigre de $Y$ $F_1 := \bigcup_{n\in J} K^1_n$ 
sur le $G_\delta$ dense de $X$, $G_1 := \bigcup_{n\in J} G^1_n$, 
o\`u $J := \{n  \in\omega ~/~ G^1_n \not= \emptyset\}$. Alors $F := F_0 \cup G_1$, $G := F_1 \cup G_0$, et 
$$f : \left\{\!\!
\begin{array}{ll} F 
& \!\!\!\rightarrow G \cr x 
& \!\!\!\mapsto \left\{\!\!\!\!\!\!
\begin{array}{ll} 
& g_0(x)~~{\rm si}~~x \in F_0 \cr & \cr
& g_1^{-1}(x)~~{\rm si}~~x \in G_1
\end{array}\right.
\end{array}\right.$$ 
r\'epondent au probl\`eme.$\hfill\square$\bigskip

 Venons-en \`a l'\'etude de la r\'eciproque du lemme de l'introduction ; dans 
ce lemme, il est question de fonctions continues et 
ouvertes, donc non n\'ecessairement injectives. Le 
contre-exemple du th\'eor\`eme 1.4 n'est donc pas un 
obstacle \`a la r\'eciproque. De fait, l'application d\'efinie par 
$$\phi : \left\{
\begin{array}{ll} \{K \in 
{\cal K}(2^\omega)\setminus \{\emptyset \} ~/~ K~\mbox{est~parfait} \} 
& \!\!\!\rightarrow P_\infty \cr K 
& \!\!\!\mapsto\mbox{max}(K)
\end{array}\right.$$ 
est surjective continue et ouverte et son graphe est contenu dans 
$$\{(K,x) \in {\cal K}(2^\omega) \setminus \{ \emptyset \} \times 2^\omega~/~x \in K \}.$$ 
Cependant, $\phi$ n'est injective sur aucun $G_\delta$ dense. Plus g\'en\'eralement, on a la
 
\begin{prop} Soient $X$ un espace polonais, 
$Y$ un espace s\'epar\'e, $f: X \rightarrow Y$ continue et 
ouverte. Si $f$ n'est pas injective, $f$ n'est injective 
sur aucun $G_\delta$ dense de $X$.\end{prop}

\noindent\bf D\'emonstration.\rm\ Posons $E := \{(x,x') \in X^2 ~/~x\not= x'~
{\rm et}~f(x) = f(x')\}$. Alors $E$ est $G_\delta$, et 
l.p.o., car si $U$ et $V$ sont ouverts dans $X$, $\Pi[E \cap 
(U\times V)] = U\cap \{x \in X ~/~ f(x) \in f[V \setminus 
\{x\}]\}$, comme on le v\'erifie imm\'ediatement. Si $f(x_0) 
\in f[V \setminus \{x_0\}]$, soit $x_1 \in V \setminus 
\{x_0\}$ tel que $f(x_0) = f(x_1)$. On trouve des ouverts 
disjoints $W_i$ de $X$ tels que $x_i \in W_i$, 
$W_1 \subseteq V$. Si $x\in f^{-1}(f[W_1]) \cap W_0$ (qui 
est un ouvert contenant $x_0$), $f(x) = f(y)$, o\`u $y \in 
W_1 \subseteq V$, et $x\not= y$ car $W_0 \cap W_1 = \emptyset$. Si $E$ est non 
vide, on applique le lemme 4.4 de [Le1] pour voir que $E$ 
rencontre tout carr\'e $G_{\delta}$ dense ; d'o\`u le r\'esultat.$\hfill\square$\bigskip

Le lemme 1.3 montre qu'on ne peut pas avoir l'image ``grosse" en g\'en\'eral, dans une uniformisation du type von Neumann, m\^eme pour un ferm\'e l.p.o. non vide d'un produit d'espaces polonais parfaits de dimension 0. Cependant, l'uniformisation souhait\'ee a lieu dans 
au moins un sens (on l'a vu avant 1.11 dans le cas 
particulier). Le th\'eor\`eme 1.13 qui suit est le r\'esultat essentiel 
de ce chapitre. Il est \`a noter que malgr\'e des hypoth\`eses 
compl\`etement sym\'etriques, la conclusion ne l'est pas (cf 
le lemme 1.3).

\begin{lem} Soient $X$ et $Y$ des espaces 
polonais, $A$ un $G_\delta$ l.p.o. non vide de $X \times Y$ 
de projections $X$ et $Y$, et $F$ (resp. $G$) un $G_\delta$ 
dense de $X$ (resp. $Y$). Alors les projections de 
$A \cap (F \times G)$ sont co-maigres.\end{lem}

\noindent\bf D\'emonstration.\rm\ Si $\Pi_X[A \cap (F \times G)]$ n'est pas 
co-maigre, soit $U$ un ouvert non vide tel que 
$$M := \Pi_X[A \cap (F \times G)] \cap U$$ 
soit maigre ; $A\cap (U \times Y)$ est $G_\delta$ l.p.o. non vide de 
$U \times Y$, donc par le lemme 4.4 de [Le1] rencontre $(U \cap 
F\setminus M) \times G$ en un point $(x,y)$ qui v\'erifie 
$x \in M\setminus M$.$\hfill\square$\bigskip

\begin{thm} Soient $X$ et $Y$ des espaces 
polonais parfaits de dimension 0, $A$ un $G_\delta$ l.p.o. 
non vide de $X \times Y$. Alors il existe des ensembles 
presque-ouverts non vides $F$ et $G$, l'un contenu dans $X$ 
et l'autre dans $Y$, et une surjection continue ouverte de 
$F$ sur $G$ dont le graphe est contenu dans $A$ ou dans
$\{(y,x) \in Y \times X~ /~ (x,y) \in A\}$ selon le 
cas.\end{thm}

\noindent\bf D\'emonstration.\bigskip

\noindent Premier cas.\rm\ Dans tout rectangle 
ouvert non vide $U \times V$ tel que les projections de 
$A \cap (U\times V)$ soient denses dans $U$ et $V$, on 
trouve un sous-rectangle $U' \times V'$ de $U \times V$ 
ayant ces propri\'et\'es et tel que pour toute partie rare $R$ 
de $V'$, $\Pi_X[A \cap (U' \times R)]$ n'est pas dense 
dans $U'$.\bigskip

 Soient $(O_n)$ une suite d'ouverts de $X \times Y$ telle 
que $A = \bigcap_{n\in\omega} O_n$, et $U_\emptyset$ et 
$V_\emptyset$ fournis par la propri\'et\'e pr\'ec\'edente 
appliqu\'ee aux projections de $A$.

\vfill\eject

 On construit alors des 
suites d'ouverts non vides $(U_s)_{s\in\omega^{<\omega}}$ 
et $(V_s)_{s\in\omega^{<\omega}}$ v\'erifiant 
$$\begin{array}{ll} 
& (i)~~~\bigcup_{n\in\omega} U_{s^\frown n}~\mbox{est~dense~dans}~U_s \cr 
& (ii)~~U_s \times V_s \subseteq O_{\vert s\vert -1}~\mbox{si}~s \not= \emptyset \cr 
& (iii)~\delta(U_s),~\delta(V_s) < \vert s\vert ^{-1}~\mbox{si}~s \not= \emptyset \cr 
& (iv)~~U_{s^\frown n} \cap U_{s^\frown m} = \emptyset~\mbox{si}~n\not= m \cr 
& (v)~~~\overline{U_{s^\frown n}} \subseteq U_s,~\overline{V_{s^\frown n}} \subseteq V_s \cr 
& (vi)~~\mbox{Les~projections~de}~A \cap (U_s \times V_s)~
\mbox{sont~denses~dans}~U_s~\mbox{et}~V_s \cr 
& (vii)~\mbox{Si}~R~\mbox{est~rare~dans}~V_s, \Pi_X[A \cap 
(U_s \times R)]~\mbox{n'est~pas~dense~dans}~U_s
\end{array}$$
Admettons la construction effectu\'ee pour $\vert s\vert  \leq p$. 
Partitionnons $\Pi_X[A \cap (U_s \times V_s)]$ 
(\'eventuellement priv\'e d'un point) en une infinit\'e 
d'ouverts-ferm\'es non vides, disons $(Z_n)$, et soit 
$$T_n := \Pi_Y[A \cap (Z_n \times V_s)]$$ 
(c'est un ouvert non vide de $V_s$). Soit $H_n$ l'ensemble des parties de 
$\boraone\lceil Z_n \setminus \{\emptyset \} \times 
\boraone\lceil Z_n \setminus \{\emptyset \}$ telles que si 
$(U, V)$ et $(U', V')$ sont distincts dans $P$, $U$ et $U'$ 
soient disjoints et $U \times V \subseteq O_{\vert s\vert }$, 
$\delta(U)$, $\delta(V) < (\vert s\vert +1)^{-1}$, $\overline{U} 
\subseteq Z_n$, $\overline{V} \subseteq T_n$, les 
projections de $A \cap (U \times V)$ soient denses dans 
$U$ et $V$, et pour toute partie rare $R$ de $V$, 
$\Pi_X[A \cap (U \times R)]$ ne soit pas dense dans $U$.\bigskip

 Alors $H_n$ n'est pas vide, puisqu'il contient le vide, 
et est ordonn\'e de fa\c con inductive par l'inclu-sion, donc 
par le lemme de Zorn a un \'el\'ement maximal $P_n$. Alors si 
$P_n := \{(U_m,V_m)~/~ m \in I_n\}$, l'ouvert $u 
:= \bigcup_{m\in I_n} U_m$ est dense dans $Z_n$, sinon soit 
$U$ un ouvert non vide de $Z_n$ disjoint de $u$, et 
$(x, y)$ dans $A \cap (U \times T_n)$ ; on trouve des 
ouverts-ferm\'es $U''$ et $V''$, de diam\`etre au plus 
$(\vert s\vert +1)^{-1}$, tels que $(x,y) \in U'' \times V'' 
\subseteq O_{\vert s\vert } \cap (U \times T_n)$. On applique 
alors la propri\'et\'e du premier cas aux projections de 
$A \cap (U'' \times V'')$ pour obtenir la contradiction 
cherch\'ee avec la maximalit\'e de $P$.\bigskip

 On obtient maintenant les $U_{s^\frown n}$ en 
renum\'erotant la suite des ouverts $U$ pour lesquels on 
trouve $V$ et $n$ tels que $(U, V) \in P_n$, et $V_{s^\frown n}$ 
est le $V$ correspondant. La construction est donc 
possible.\bigskip

 Soit $F' := \bigcap_{n\in\omega} \bigcup_{s\in \omega^n} 
U_s$ ; 
alors $F'$ est $G_\delta$ dense de $U_\emptyset$, donc 
presque-ouvert. Si $x$ est dans $F'$, on d\'efinit $f(x)$ de 
la mani\`ere habituelle, et la fonction $f$ est continue et 
uniformise partiellement $A$. De plus, si un ouvert non 
vide de $F'$ avait une image maigre, il contiendrait un 
ouvert non vide d'image rare par le th\'eor\`eme de Baire et 
la continuit\'e de $f$. Mais ceci est exclus \`a cause de la 
condition (vii), puisque $f[F' \cap U_s] \subseteq V_s$.\bigskip

 Soit $(U_n)$ une base de la topologie de $F'$, $(V_n)$ une 
suite d'ouverts de $Y$, et $M_n$ une suite de $F_\sigma$ 
maigres de $Y$ tels que $f[F' \cap U_n]\Delta V_n \subseteq M_n$. 
Posons $F = F' \setminus (\bigcup_{n\in\omega} f^{-1}(M_n))$ ; $F$ est $G_\delta$ 
dense de $U_\emptyset$ comme $F'$ par ce qui pr\'ec\`ede. Il 
est maintenant clair qu'on peut poser 
$$G := (\bigcup_{n\in\omega} V_n) \setminus (\bigcup_{n\in\omega} M_n).$$
\bf Second cas.\rm\ Il existe un rectangle ouvert non 
vide $U_\emptyset \times V_\emptyset$ tel que les 
projections de $A \cap (U_\emptyset \times V_\emptyset)$ 
soient denses dans $U_\emptyset$ et $V_\emptyset$, tel que 
pour tout sous-rectangle $U' \times V'$ de 
$U_\emptyset \times V_\emptyset$ ayant ces propri\'et\'es, on 
trouve une partie rare $R$ de $V'$ telle que 
$\Pi_X[A\cap (U' \times R)]$ soit dense dans $U'$.\bigskip

 Soient $n_0 \in \omega$, $\varepsilon > 0$, et $U'$ et $V'$ des 
ouverts ayant ces propri\'et\'es ; on montre qu'il \noindent{existe} des suites d'ouverts non vides $(U_n)$ 
et $(V_n)$ telles que 
$$\begin{array}{ll} 
& (i)~~~\bigcup_{n\in \omega} U_n~
\mbox{est~dense~dans~} U',~\bigcup_{n\in \omega} V_n~
\mbox{est~dense~dans~} V' \cr  
& (ii)~~U_n\times V_n \subseteq 
O_{n_0} \cr 
& (iii)~\delta(U_n),~\delta(V_n)< \varepsilon \cr 
& (iv)~~V_n\cap V_m = \emptyset ~\mbox{si}~n\not= m \cr 
& (v)~~~\overline{U_n} \subseteq U',~
\overline{V_n} \subseteq V' \cr 
& (vi)~~\mbox{Les~projections~de~} A\cap (U_n\times V_n)~\mbox{sont~
denses~dans~} U_n ~\mbox{et~} V_n
\end{array}$$ 
 Soit $(x_n)$ une suite dense de 
$\Pi_X[A\cap (U'\times V')]$ (donc de U'). On va commencer 
par construire des suites d'ouverts-ferm\'es non vides 
$(Z_n)$ et $(T_n)$ v\'erifiant $A\cap (Z_n\times T_n) \not= \emptyset$, 
$$Z_n\times T_n \subseteq O_{n_0}\cap 
[(U'\cap {\cal B}(x_n,2^{-n}))\times V'\setminus 
(\bigcup_{p<n} T_p \cup R)]\mbox{,}$$ 
de diam\`etre au plus $\varepsilon$.\bigskip

 Admettons avoir trouv\'e $(Z_p)_{p<n}$ et $(T_p)_{p<n}$ ayant 
ces propri\'et\'es. Comme $R$ est rare et non vide, 
$\bigcup_{p<n} T_p \subset_{\not=} V'\setminus R$, et les 
projections de $A\cap [U'\times (V'\setminus 
(\bigcup_{p<n} T_p))]$ sont co-maigres dans $U'$ et 
$V'\setminus (\bigcup_{p<n} T_p)$. Par le lemme pr\'ec\'edent, 
celles de $A_n := A\cap [U'\times (V'\setminus 
(\bigcup_{p<n} T_p \cup R))]$ le sont \'egalement ; la 
projection sur $X$ rencontre donc ${\cal B}(x_n,2^{-n})$ en 
$z_n$ ; soit $y_n$ tel que $(z_n,y_n)\in A_n$, et $Z_n$, $T_n$ 
des ouverts-ferm\'es de diam\`etre au plus $\varepsilon$ tels que l'on ait 
$$(z_n,y_n)\in Z_n\times T_n \subseteq O_{n_0}\cap 
[(U'\cap {\cal B}(x_n,2^{-n}))\times V'\setminus 
(\bigcup_{p<n} T_p \cup R)].$$ 
La construction est donc possible.\bigskip

 Si $\bigcup_{n\in \omega} \Pi_Y[A\cap (Z_n\times T_n)]$ est 
dense dans $V'$, on d\'efinit 
$V_n := \Pi_Y[A\cap (Z_n\times T_n)]$ ainsi que 
$$U_n := \Pi_X[A\cap (Z_n\times V_n)]\mbox{,}$$ 
et les conditions sont v\'erifi\'ees, par densit\'e de la suite $(z_n)$ dans $U'$. Sinon, 
la construction pr\'ec\'edente montre que les autres conditions 
sont r\'ealis\'ees. On pose 
$$Y' := \Pi_Y[A\cap (U'\times V')] 
\setminus \overline{\bigcup_{n\in \omega} 
\Pi_Y[A\cap (Z_n\times T_n)]}.$$ 
Si $(x,y) \in A\cap (U'\times Y')$, soit $Z_{x,y}\times T_{x,y}$ un rectangle ouvert de diam\`etre 
au plus $\varepsilon$ tel que l'on ait les inclusions $(x,y)\in 
Z_{x,y}\times T_{x,y} \subseteq \overline{Z_{x,y}\times 
T_{x,y}} \subseteq O_{n_0}\cap (U'\times Y')$ et tel que 
les projections de $A\cap (Z_{x,y}\times T_{x,y})$ soient 
$Z_{x,y}$ et $T_{x,y}$. On a $Y' = \bigcup_{x,y} T_{x,y} = 
\bigcup_{n\in \omega} T_{x_n,y_n}$. R\'eduisons la suite 
$(T_{x_n,y_n})$ en $(T'_n)$, et posons $Z'_n := \Pi_X[A
\cap (Z_{x_n,y_n}\times T'_n)]$. Alors en ne gardant que 
les $T'_n$ non vides, les $Z'_n$ correspondants, et les 
$Z_n$ et $T_n$, on a ce qu'on veut.

\vfill\eject

 On recopie alors quasiment la construction du th\'eor\`eme 5.2 de [Le1] : 
soit $\phi_0$ de $\omega$ dans $\{ \emptyset \}$ et, si $n>0$, 
$\phi_n$ une bijection de $\omega$ sur $\omega^n$. On construit 
une suite $(U_s)_{s\in \omega^{<\omega}}$ d'ouverts non vides 
de $X$, et une suite $(V_s)_{s\in \omega^{<\omega}}$ d'ouverts non vides 
de $Y$ v\'erifiant 
$$\begin{array}{ll} 
& (i)~~~\bigcup_{n\in \omega} 
\overline{U_{s^\frown n}}~
\mbox{est~dense~dans~} U_s,~\bigcup_{n\in \omega} 
\overline{V_{s^\frown n}}~
\mbox{est~dense~dans~} V_s \cr  
& (ii)~~U_s\times V_s \subseteq 
O_{\vert s\vert -1}~\mbox{si~} s\not= \emptyset  \cr 
& (iii)~\delta(U_s),~\delta(V_s)< \vert s\vert ^{-1} ~\mbox{si~} s\not= \emptyset\cr 
& (iv)~~V_{s^\frown n}\cap V_{s^\frown m} = \emptyset 
~\mbox{si}~n\not= m \cr 
& (v)~~~(\bigcup_{k\in \omega} U_{\phi_n(k)})\cap 
\bigcup_{q+p<n} [U_{\phi_q(p)}\setminus 
(\bigcup_{l\in \omega} U_{\phi_q(p)^\frown l})] = \emptyset\cr 
& (vi)~~\mbox{Les~projections~de~} A\cap (U_s\times V_s)~\mbox{sont~
denses~dans~} U_s ~\mbox{et~} V_s
\end{array}$$
Admettons avoir construit les suites $(U_s)_{\vert s\vert  
\leq n}$ et $(U_s)_{\vert s\vert  \leq n}$, 
$(U_{\phi_n(p)^\frown k})_{p<m,k\in \omega}$, 
$(V_{\phi_n(p)^\frown k})_{p<m,k\in \omega}$ v\'erifiant (i)-(vi).\bigskip

 On construit, si ce n'est d\'ej\`a fait, 
$(U_{\phi_n(m)^\frown k})_{k\in \omega}$ et 
$(V_{\phi_n(m)^\frown k})_{k\in \omega}$ en appliquant ce 
qui pr\'ec\`ede \`a $\varepsilon = (n+1)^{-1}$, $V' := V_{\phi_n(m)}$, 
$U' := U_{\phi_n(m)}\setminus \overline{\bigcup_{q+p<n} 
[U_{\phi_q(p)}\setminus (\bigcup_{l\in \omega} 
U_{\phi_q(p)^\frown l})]}$.\bigskip

 Les conditions demand\'ees sont v\'erifi\'ees, la densit\'e des 
projections de $A\cap (U'\times V')$ dans $U'$ et $V'$, 
donc dans $U_{\phi_n(m)}$ et $V_{\phi_n(m)}$, ne posant pas 
de probl\`eme \`a cause du lemme pr\'ec\'edent. On d\'efinit alors 
$F$ et $G$ comme dans le th\'eor\`eme 5.2 de [Le1] ; ce sont des 
$G_{\delta}$ denses de $V_{\emptyset}$ et $U_{\emptyset}$, 
donc des presque-ouverts non vides. On conclut alors comme 
dans le th\'eor\`eme 5.2 de [Le1].$\hfill\square$\bigskip

 Sous les hypoth\`eses du th\'eor\`eme 1.13, il est faux en 
g\'en\'eral que $A$ est uniformisable sur un $G_{\delta}$ dense 
de sa projection par une application continue ouverte sur son 
image (bien que $A$ soit uniformisable par une application 
continue). En effet, on prend l'exemple du lemme 1.3, ce 
qui fournit $A_0 \subseteq N_{(0)}\times \omega^\omega$ ; si $A_1$ 
est le graphe d'un hom\'eomorphisme de $\omega^\omegaÊ\setminus N_{(0)}$ 
sur $\omega^\omega$, $A_0 \cup A_1$ est ferm\'e l.p.o. de 
$\omega^\omega \times \omega^\omega$, de projections $\omega^\omega$. Raisonnons 
par l'absurde : il existe un $G_{\delta}$ dense $G$ de $\omega^\omega$ 
et une application $f$ continue sur $G$ et ouverte sur 
$f[G]$, qui uniformise $A$ partiellement. Alors 
$f[G\cap N_{(0)}]$ est $G_{\delta}$ rare, par construction 
de $A_0$.\bigskip

 Mais comme $f[G\setminus N_{(0)}]$ est $G_{\delta}$ 
dense de $\omega^\omega$, $f[G]$ aussi, donc $f[G\cap N_{(0)}]$ est 
ouvert non vide et rare de $f[G]$, ce qui est absurde. 
Cependant, l'uniformisation a lieu si $A$ est le graphe d'une 
surjection continue ouverte de $Y$ dans $X$ :

\begin{prop} Soient $X$ et $Y$ des espaces 
m\'etrisables s\'eparables de dimension 0, $Y$ \'etant complet, et 
$f:Y\rightarrow X$ une surjection continue ouverte ; alors 
il existe un hom\'eomorphisme $g$ de $X$ sur son image tel que
 $f \circ g = \mbox{Id}_X$.\end{prop}

\noindent\bf D\'emonstration.\rm\ On construit des suites d'ouverts-ferm\'es 
$(U_s)_{s\in \omega^{<\omega}}$ et 
$(V_s)_{s\in \omega^{<\omega}}$ v\'erifiant 
$$\begin{array}{ll} 
& (i)~~~\bigcup_{n\in \omega} 
U_{s^\frown n} = U_s \cr  
& (ii)~~U_s = f[V_s] \cr 
& (iii)~\delta(U_s),~\delta(V_s)< \vert s\vert ^{-1} ~\mbox{si~} 
s\not= \emptyset \cr 
& (iv)~~U_{s^\frown n}\cap U_{s^\frown m} = 
V_{s^\frown n}\cap V_{s^\frown m} = \emptyset ~\mbox{si}~n
\not= m \cr 
& (v)~~~V_{s^\frown n} \subseteq V_s
\end{array}$$

\vfill\eject

 On pose $U_{\emptyset} := X$, $V_{\emptyset} := Y$ ; si 
on a $U_s$ et $V_s$, on partitionne $U_s$ en une suite 
$(U_n)$ d'ouverts-ferm\'es de diam\`etre au plus 
$(\vert s\vert +1)^{-1}$ ; on partitionne ensuite l'ouvert-ferm\'e 
$V_s\cap f^{-1}(U_n)$ en une suite $(V_m^n)_m$ 
d'ouverts-ferm\'es de diam\`etre au plus $(\vert s\vert +1)^{-1}$. La suite 
double d'ouverts $(f[V_m^n])$ est ensuite r\'eduite en $(W_m^n)$, 
et on renum\'erote les $W_m^n$, l'ouvert-ferm\'e correspondant 
 dans $Y$ \'etant $V_m^n\cap f^{-1}(W_m^n)$.\bigskip
 
 Cette construction \'etant faite, elle d\'efinit de la mani\`ere 
habituelle une application continue 
$$g:X\rightarrow Y.$$
De plus, si $\alpha$ est tel que 
$x\in \bigcap_{n\in \omega} U_{\alpha \lceil n}$, on a 
$f(g(x)) \in \bigcap_{n\in \omega} 
f[V_{\alpha \lceil n}] = \bigcap_{n\in \omega} U_{\alpha 
\lceil n} = \{ x \}$, d'o\`u le fait que $f\circ g =\mbox{Id}_X$ et 
l'injectivit\'e de $g$. Enfin, par la condition (ii), on a 
$g[U_s] = g[X]\cap V_s$, donc $g$ est ouverte sur son image.$\hfill\square$\bigskip
 
 Bien s\^ur, ce r\'esultat est faux si on enl\`eve la condition de 
dimension sur $X$.

\begin{cor} Soient $X$ et $Y$ des espaces 
polonais parfaits de dimension 0, $A$ un $G_{\delta}$ l.p.o. 
non vide de $X\times Y$. Alors $A$ est uniformisable sur un 
presque-ouvert non vide de $X$ par une application continue 
et ouverte sur son image.\end{cor}

\noindent\bf D\'emonstration.\rm\ On applique le th\'eor\`eme 1.13 : on a le 
r\'esultat tout de suite ou alors on applique la proposition 
pr\'ec\'edente.$\hfill\square$\bigskip

 A cause du lemme 1.3, l'image de l'application fournie par ce corollaire est en g\'en\'eral rare.

\section{$\!\!\!\!\!\!$ Applications aux classes de Wadge potentielles.}\indent
 
 On va maintenant appliquer les r\'esultats d'uniformisation 
aux classes de Wadge {potentielles} : on obtient pour commencer des 
caract\'erisations des ensembles potentiellement diff\'erence
 d'ouverts parmi les ensembles potentiellement 
$\borathree$ et $\bormthree$ (donc par exemple parmi les bor\'eliens \`a coupes d\'enombrables).

\begin{lem} Soit $X$ un espace polonais r\'ecursivement 
pr\'esent\'e. Alors il existe un $G_{\delta}$ dense $\Ana$ de 
$X$ sur lequel la topologie de $X$ co\"\i ncide avec la topologie 
$\it\Delta$ engendr\'ee par les $\Borel$.\end{lem}

\noindent\bf D\'emonstration.\rm\ Il est prouv\'e dans [Ke] que si $A$ 
est $\Borel$, il existe un $\Borel$ ouvert $U$ et un $\Borel$ 
\`a coupes ferm\'ees rares $F$ de $\omega \times X$ tels que 
$A \Delta U \subseteq \bigcup_{p \in \omega} F_p$. Cette 
propri\'et\'e s'\'ecrit de mani\`ere $\Ca$ ; on peut donc associer, 
\`a tout entier $n$ codant un $\Borel$, un entier $(f(n))_0$ 
codant un ouvert $\Borel$, et un entier $(f(n))_1$ codant 
un $\Borel$ \`a coupes ferm\'ees rares de $\omega \times X$, en 
gardant l'inclusion ci-dessus, et ce par une fonction 
$\Ca$-r\'ecursive partielle $f$. D\'esignons par $W^X$ un 
ensemble $\Ca \subseteq \omega$ de codes pour les $\Borel$ 
de $X$ et par $C^X \subseteq \omega \times X$ un ensemble 
$\Ca$ dont les sections aux points de $W^X$ d\'ecrivent les 
$\Borel$ de $X$ (cf [Lo2]). 

\vfill\eject 

 On pose $x \in G \Leftrightarrow \forall~n~\forall~p~~n \notin W^X$ 
ou $((f(n))_1,~p,~x) \notin C^{\omega \times X}$. Alors $G$ 
r\'epond au probl\`eme, car $\check G$ est r\'eunion d\'enombrable de 
ferm\'es rares, et si $n$ code $A\in \Borel$, on a 
$A\cap G = C^X_{(f(n))_0}\cap G$, qui est ouvert dans $G$.$\hfill\square$\bigskip
 
 Dans la suite, si $f_s$ est une fonction partielle de $X$ 
dans $Y$ ou de $Y$ dans $X$, on notera $G(f_s)$ la partie de 
$X\times Y$ \'egale au graphe de $f_s$ si $f_s$ va de $X$ 
dans $Y$, et \`a $\mbox{Gr}(f_s)^* := \{ (x,y)~/~x=f_s(y) \}$ sinon.
 
\begin{lem} Soient $X$ un espace polonais parfait 
r\'ecursivement pr\'esent\'e, $G$ un universel pour les $\Ca$ de 
$X$, $\varphi$ et $\psi$ des $\Ca$-op\'erateurs 
 monotones sur $X$, et $\Phi^\xi$ les op\'erateurs sur $X$ 
d\'efinis par $\Phi^0(A) = A$, et 
$\Phi^\xi(A) = \varphi (\bigcup_{\eta < \xi} \Phi^\eta (A))$ 
si $\xi$ est impair, $\Phi^\xi(A) = \psi (\bigcup_{\eta < \xi} \Phi^\eta (A))$ 
si $\xi > 0$ est pair. Alors si $R(\alpha ,\beta ,x) 
\Leftrightarrow \alpha \in WO$ et 
$x \in \Phi^{\vert\alpha\vert }(G_{\beta})$, la relation $R$ 
est $\Ca$.\end{lem}

\noindent\bf D\'emonstration.\rm\ Ce lemme est d\'emontr\'e dans [Lo1] si 
$\varphi = \psi = \Phi$. La d\'emonstration ici est analogue. Pr\'ecisons-en les diff\'erences.\bigskip

 On sait qu'il existe une fonction $\Borel$-r\'ecursive 
$g:\omega^\omega \times \omega \rightarrow \omega^\omega$ telle que 
si $\alpha \in WO$, $g(\alpha,n)$ est dans $WO$ et code 
l'odre $\leq_{\alpha}$ restreint aux pr\'ed\'ecesseurs de $n$ 
pour $\leq_{\alpha}$. On d\'efinit par r\'ecurrence une 
fonction $\Borel$-r\'ecursive $\varepsilon : \omega^\omega \rightarrow 
\omega$ telle que si $\alpha \in WO$, on ait l'\'egalit\'e 
$$\varepsilon (\alpha ) = \left \{\!\!\!\!\!\!
\begin{array}{ll}  
& 0~\mbox{si}~\vert\alpha\vert ~\mbox{est~pair,} \cr 
& 1~\mbox{si}~\vert\alpha\vert ~\mbox{est~impair.}
\end{array}\right.$$ 
On d\'efinit alors un $\Ca$-op\'erateur monotone $\Psi$ par 
$$(n,\alpha ,\beta ,x) \in \Psi(P)\ \Leftrightarrow\left\{\!\!\!\!\!\!
\begin{array}{ll}
& (n=0\mbox{ et }\alpha \in WO\mbox{ et }\forall~p~~(0,g(\alpha ,p),\beta ,x) \in P)\mbox{ ou}\cr
& (n=1$ et $\alpha \in WO\mbox{ et }\exists ~q~/~q \leq_{\alpha} q\mbox{ et }\forall~p~(0,g(\alpha ,p),\beta ,x) \in P\mbox{ et}\cr
& [(\varepsilon (\alpha )=0\mbox{ et }x \in \psi(\{y \in X~/~\exists~q~P(1,g(\alpha ,q),\beta ,y) \}))\mbox{ ou}\cr
& (\varepsilon (\alpha )=1\mbox{ et }x \in \varphi(\{y \in X~/~\exists~q~P(1,g(\alpha ,q),\beta ,y) \}))])
\mbox{ ou}\cr
& P(n,\alpha ,\beta ,x)
\end{array}\right.$$
On montre alors comme dans [Lo1] que 
$$\begin{array}{ll}
(0,\alpha ,\beta ,x) \in \Psi^{\xi}(P_0) 
& \!\!\!\!\Leftrightarrow\alpha \in WO\mbox{ et }\vert\alpha\vert  \leq \xi\mbox{, et}\cr
(1,\alpha ,\beta ,x) \in \Psi^{\xi}(P_0)
& \!\!\!\!\Leftrightarrow\alpha \in WO\mbox{ et }\vert\alpha\vert  \leq \xi
\mbox{ et }x \in \Phi^{\vert\alpha\vert }(G_{\beta})\mbox{, o\`u}\cr 
\ \ \ \ \ \ \ \ \ \ \ \ P_0(n,\alpha ,\beta ,x)
& \!\!\!\!\Leftrightarrow\left\{\!\!\!\!\!\!\! 
\begin{array}{ll}
& (n=0~\mbox{et}~\alpha \in WO~\mbox{et}~\forall~q~~q \not \leq _{\alpha} q)~\mbox{ou}\cr 
& (n=1~\mbox{et}~\alpha \in WO~\mbox{et}~\forall~q~~q \not\leq _{\alpha} q~\mbox{et}~G(\beta ,x)).
\end{array}\right.
\end{array}$$
Comme l'op\'erateur $\Psi$ et $P_0$ sont $\Ca$, on a donc 
que $\Psi^{\infty}(P_0)$ est $\Ca$, et aussi que $(1,\alpha ,\beta ,x)$ est dans $\Psi^{\infty}(P_0)$ ssi 
($\alpha \in WO$ et $x \in \Phi^{\vert\alpha\vert }(G_{\beta })$). On a donc le 
r\'esultat, puique $R = \Psi^{\infty}(P_0)_1$.$\hfill\square$\bigskip

 Soit $\xi$ un ordinal d\'enombrable non nul. On d\'efinit 
$f:\omega^{<\omega} \rightarrow \{ -1 \} \cup (\xi +1)$, par 
r\'ecurrence sur $\vert s\vert $, comme suit : $f(\emptyset) = \xi$ et 
$$f(s^\frown n) = \left\{\!\!\!\!\!\!
\begin{array}{ll} 
& \bullet~-1~\mbox{si}~f(s)\leq 0\mbox{,}\cr 
& \bullet~\theta~\mbox{si}~f(s)=\theta +1\mbox{,}\cr
& \bullet~\mbox{un~ordinal~impair~de}~f(s)~\mbox{tel~que~la~suite}~
(f(s^\frown n))_n~\mbox{soit~co-finale~dans} \cr & f(s)~\mbox{et~
strictement~croissante~si}~f(s)~\mbox{est~limite~non~nul.}
\end{array}\right.$$
On d\'efinit alors des arbres : 
$T_{\xi} := \{ sÊ\in \omega^{<\omega}~/~f(s) \not= -1 \}$ et 
$T'_{\xi} := \{ sÊ\in T_{\xi}~/~f(s) \not= 0 \}$.

\vfill\eject

\begin{thm} Soient $X$ et $Y$ des espaces 
polonais, $A$ un bor\'elien $\mbox{pot}(\borathree)\cap \mbox{pot}(\bormthree)$ de 
$X \times Y$, et $\xi$ un ordinal d\'enombrable.\smallskip

\noindent (a) Si $\xi$ est pair non nul, $A$ est non-$\mbox{pot}(D_{\xi}(\boraone))$ 
si et seulement s'il existe des espaces polonais parfaits 
$Z$ et $T$ de dimension 0, des ouverts-ferm\'es non vides 
$A_s$ et $B_s$ (l'un dans $Z$ et l'autre dans $T$, pour 
$s$ dans $T_{\xi}$), des surjections continues ouvertes 
$f_s$ de $A_s$ sur $B_s$, et des injections continues $u$ 
et $v$ tels que si 
$B_p := \bigcup_{s \in T_{\xi}~/~
\vert s\vert ~\mbox{paire}} G(f_s)$ et ${B_i := \bigcup_{s \in T_{\xi}~/~
\vert s\vert ~\mbox{impaire}} G(f_s)}$, on ait $\overline{B_p} = B_p \cup B_i$, 
$B_p \subseteq (u\times v)^{-1}(A)$, $B_i \subseteq (u\times v)^{-1}(\check A)$, 
et ${G(f_s) = \overline{\bigcup_{n \in \omega} G(f_{s^\frown n})} \setminus 
(\bigcup_{n \in \omega} G(f_{s^\frown n}))}$ si $s \in T'_{\xi}$.\medskip

\noindent (b) Si $\xi$ est impair, $A$ est non-$\mbox{pot}(\check D_{\xi}(\boraone))$ si et 
seulement s'il existe des espaces polonais $Z$ 
et $T$ parfaits de dimension 0, des ouverts-ferm\'es non vides $A_s$ 
et $B_s$ (l'un dans $Z$ et l'autre dans $T$, pour $s$ dans 
$T_{\xi}$), des surjections continues ouvertes 
$f_s$ de $A_s$ sur $B_s$, et des injections continues $u$ 
et $v$ tels que si $B_p := \bigcup_{s \in T_{\xi}~/~
\vert s\vert ~\mbox{paire}} G(f_s)$ et $B_i := \bigcup_{s \in T_{\xi}~/~
\vert s\vert ~\mbox{impaire}} G(f_s)$, on ait $\overline{B_i} = B_p \cup B_i$, 
$B_i \subseteq (u\times v)^{-1}(A)$, $B_p \subseteq (u\times v)^{-1}(\check A)$, 
et ${G(f_s) = \overline{\bigcup_{n \in \omega} G(f_{s^\frown n})} \setminus 
(\bigcup_{n \in \omega} G(f_{s^\frown n}))}$ si $s \in T'_{\xi}$.\end{thm}

\noindent\bf D\'emonstration.\rm\ Montrons (a), la preuve de (b) \'etant analogue. 
Supposons que $A$ est $\mbox{pot}(D_{\xi}(\boraone))$, en raisonnant 
par l'absurde ; alors $B_p = \overline{B_p} \cap (u\times v)^{-1}(A)$, 
l'est aussi, cette classe \'etant stable par intersection avec les ferm\'es, 
et on trouve un $G_{\delta}$ dense $F$ (resp. $G$) de $Z$ (resp. $T$) tels 
que $B_p\cap (F\times G)$ soit $D_{\xi}(\boraone)$ dans 
$F\times G$. On trouve donc une suite croissante d'ouverts 
de $F\times G$, disons $(U_{\eta})_{\eta<\xi}$, telle que 
$B_p\cap (F\times G) = \bigcup_{\eta<\xi,\eta~\mbox{impair}} 
U_{\eta} \setminus (\bigcup_{\theta<\eta} U_{\theta})$.\bigskip

 Montrons que si $\eta \leq \xi$, $s \in T_{\xi}$ et 
$f(s) = \eta$, alors $G(f_s)\cap (F\times G) \subseteq 
\check U_{\eta}$ si $\eta < \xi$, et $G(f_s)\cap 
(F\times G) \subseteq ~ ^c(\bigcup_{\theta<\eta} U_{\theta})$ si 
$\eta = \xi$. On aura la contradiction cherch\'ee avec 
$s = \emptyset$ et $\eta = \xi$. \bigskip

 On proc\`ede par r\'ecurrence sur $\eta$. Remarquons que si 
$s$ est dans $T_{\xi}$, $\vert s\vert $ est paire si et 
seulement si $f(s)$ est pair. Si $\eta =0$, $\vert s\vert $ est 
paire, donc $G(f_s)\cap (F\times G) \subseteq B_p\cap (F\times 
G) \subseteq \check U_0$.\bigskip

 Admettons le r\'esultat pour $\theta<\eta$. Si $\eta$ est le 
successeur de $\theta$, par hypoth\`ese de r\'ecurrence on a 
$G(f_{s^\frown m})\cap (F\times G) \subseteq \check U_{\theta}$ 
pour tout $m$ ; d'o\`u, si on pose ${C_s := \bigcup_{n \in \omega} 
G(f_{s^\frown n})}$, 
$$C_s\cap (F\times G) \subseteq \check U_{\theta}$$ 
et $\overline{C_s\cap (F\times G)}^{F\times G}\subseteq \check U_{\theta}$. Mais comme dans la preuve du 
lemme 3.5 de [Le1], on a 
$$\overline{C_s\cap (F\times G)}^{F\times G} = \overline{C_s}\cap (F\times G)\mbox{,}$$ 
d'o\`u l'inclusion cherch\'ee si $\eta = \xi$.\bigskip

 Si $\eta<\xi$ et $\vert s\vert $ est paire, $f(s)$ est pair et 
$\theta$ est impair. On a les inclusions successives 
$$G(f_s)\cap (F\times G) 
\subseteq B_p\cap (F\times G) \subseteq 
\bigcup_{\varepsilon<\xi,\varepsilon~\mbox{impair}} U_{\varepsilon} \setminus 
(\bigcup_{\gamma<\varepsilon} U_{\gamma}) \subseteq \check U_{\theta +1}.$$
Si $\vert s\vert $ est impaire, $f(s)$ est impair et $\theta$ est 
pair. Mais si $s \in T'_{\xi}$ est de 
longueur impaire, on a 
$${G(f_s)\cap (F\times G) \subseteq (F\times G) 
\setminus B_p = (F\times G) \setminus (\bigcup_{\varepsilon<\xi} U_{\varepsilon}) \cup 
\bigcup_{\varepsilon<\xi,\varepsilon~\mbox{pair}} U_{\varepsilon} \setminus 
(\bigcup_{\gamma<\varepsilon} U_{\gamma})}\hbox{,}$$ 
ceci parce que $G(f_s)\cap B_p = \emptyset$. On a donc le r\'esultat en appliquant l'hypoth\`ese de r\'ecurrence.

\vfill\eject

 Si $\eta$ est un ordinal limite, $(f(s^{\frown }n))_n$ est 
co-finale dans $f(s)$, donc par hypoth\`ese de r\'ecurrence, on 
a $G(f_{s^{\frown }n})\cap (F\times G) \subseteq \check 
U_{f(s^{\frown }n)}$. Si $\theta_0<f(s)$, on trouve $n(\theta_0)$ 
 tel que $f(s^{\frown} n) > \theta_0$ si ${n(\theta_0) \leq n}$. 
Donc $G(f_{s^{\frown} n}) \cap (F\times G) \subseteq \check 
U_{\theta_0}$ d\`es que $n(\theta_0) \leq n$. Or  
$${G(f_s)\!\cap \!(F\times G)} \!\subseteq \!
(F\times G)\cap \overline{C_s} \setminus C_s\! = \!
\overline{(F\times G)\!\cap\! C_s}^{\!F\times G\!} \setminus C_s 
\!\subseteq\! \overline{\!\!\bigcup\limits_{n(\theta_0)\leq n}\!\! (\!F\times G\!)
\!\cap\! G(f_{s^{\frown} n})}^{F\times G}{\!\subseteq \!\check U_{\theta_0}}.$$ 
Donc $G(f_s)\cap (F\times G) \subseteq 
~ ^c(\bigcup_{\theta<\eta} U_{\theta})$. Si $\eta<\xi$, comme 
$\vert s\vert $ est paire, on a l'inclusion 
$$G(f_s)\cap (F\times G) 
\subseteq \bigcup_{\varepsilon<\xi,\varepsilon~\mbox{impair}} U_{\varepsilon} 
\setminus (\bigcup_{\gamma<\varepsilon} U_{\gamma})\mbox{,}$$ 
donc $G(f_s)\cap (F\times G) \subseteq \check U_{\eta}$.\bigskip

 Inversement, soit $A$ dans 
$\mbox{pot}(\borathree)\cap \mbox{pot}(\bormthree)\setminus \mbox{pot}(D_{\xi}(\boraone))$ dans 
$X\times Y$, avec $X$ et $Y$ polonais. N\'ecessairement $X$ et $Y$ sont non d\'enombrables, donc 
bor\'eliennement isomorphes \`a $\omega^\omega$, disons par $\varphi$ et $\psi$. Remarquons qu'on peut supposer $X$ et $Y$ parfaits. Soit alors $\beta$ tel que $\xi < \omega_1^{\beta}$, tel que $X$ et $Y$ soient r\'ecursivement en $\beta$-pr\'esent\'es, tel que 
$\varphi$ et $\psi$ soient $\Borel (\beta)$, et tel que $A$ et 
$\check A$ soient r\'eunion d'une $\Borel (\beta)$-suite de 
$G_{\delta}$ pour le produit ${\it\Delta}^{\beta}_X\times {\it\Delta}^{\beta}_Y$ (o\`u ${\it\Delta}^{\beta}_X$ est la topologie engendr\'ee par les $\Borel (\beta)$ de $X$). Posons 
$\Omega ^{\beta}_X := \{ x\in X~/~\omega _1^{(\varphi (x),
\beta)} \leq \omega_1^\beta \}$, 
$D ^{\beta}_X := \{ x\in X~/~x\notin \Borel (\beta)\}$, 
$Z_0 := \Omega ^{\beta}_X\cap D ^{\beta}_X$, 
$T_0 := \Omega ^{\beta}_Y\cap D ^{\beta}_Y$.\bigskip

 Ces espaces $Z_0$ et $T_0$, si on les munit des 
restrictions des topologies de Gandy-Harrington 
${\it\Sigma}^{\beta}_X$ (resp. ${\it\Sigma}^{\beta}_X$), 
sont polonais parfaits de dimension 0. En effet, $Z_0$ et 
$T_0$ sont $\Ana (\beta)$ comme intersection de deux 
$\Ana (\beta)$ (cf [Mo]). Ceci prouve qu'ils n'ont pas de point 
isol\'e. De plus, si $E$ est $\Ana (\beta)$ contenu dans 
$\Omega ^{\beta}_X$, $E$ est ouvert-ferm\'e dans 
$\Omega ^{\beta}_X$ pour la restriction de 
${\it\Sigma}^{\beta}_X$ : on a, si $f$ est $\Borel (\beta)$ 
telle que $x \notin E \Leftrightarrow f(x) \in WO$, 
$$x \in \Omega ^{\beta}_X \setminus E 
\Leftrightarrow \exists~\xi < \omega_1^\beta~
(f(x) \in WO\mbox{ et }\vert f(x)\vert  \leq \xi)\mbox{ et }x \in 
\Omega ^{\beta}_X.$$
On en d\'eduit que $\Omega ^{\beta}_X$ est \`a base 
d\'enombrable d'ouverts-ferm\'es, donc m\'etrisable s\'eparable ; 
comme il est ouvert d'un espace fortement 
$\alpha$-favorable, il est lui-m\^eme fortement 
$\alpha$-favorable, donc polonais (cf [Lo1] pour plus 
de d\'etails).\bigskip

 On d\'efinit $\Omega ^{\beta}_{X\times Y} := 
\{ (x,y)\in X\times Y~/~\omega _1^{(\varphi (x),\psi (y),
\beta)} \leq \omega_1^\beta \}$ ; alors on sait que 
$\Omega ^{\beta}_{X\times Y}$ rencontre tout ensemble 
$\Ana (\beta)$ non vide de $X\times Y$ (cf [Lo1]).\bigskip

 On d\'efinit ensuite par r\'ecurrence 
$$F_{\eta} := 
\left\{\!\!\!\!\!\!
\begin{array}{ll} 
& \overline{A \cap \bigcap_{\theta<\eta} 
F_{\theta}} ~\mbox{si}~\eta~\mbox{est~pair,} \cr & \cr
& \overline{\check A \cap \bigcap_{\theta<\eta} 
F_{\theta}} ~\mbox{ si}~\eta~\mbox{est~impair,}
\end{array}
\right.$$ 
l'adh\'erence \'etant prise au sens de la topologie 
${\it\Delta}^{\beta}_X\times {\it\Delta}^{\beta}_Y$. 

\vfill\eject

 Montrons que les $F_{\eta}$ sont $\Ana (\beta)$ pour 
$\eta \leq \xi$. On d\'efinit des $\Ca (\beta)$-op\'erateurs monotones 
sur $X\times Y$ par les formules :
 $\varphi(P) = P \cup\mbox{Int}(P \cup A)$ et 
 $\psi(P) = P \cup\mbox{Int}(P \cup \check A)$, o\`u l'int\'erieur 
est pris au sens de la topologie ${\it\Delta}^{\beta}_X\times {\it\Delta}^{\beta}_Y$. Il est manifeste que, avec les notations du lemme 2.2, on a 
$\check F_{\eta} = \Phi^{\eta}(\check {\overline{A}})$ si $\eta$ est d\'enombrable, par 
r\'ecurrence transfinie. Il suffit alors d'appliquer ce lemme 2.2.\bigskip

 Montrons que si $C$ et $D$ sont co-d\'enombrables, alors 
$F_{\xi} \cap (C\times D) \not= \emptyset$. Pour ce faire, posons 
$A' := A \cap (C\times D)$. Par la remarque 2.1 de [Le1], $A'$ est 
non-$\mbox{pot}(D_{\xi}(\boraone))$. Soit $\gamma$ tel que $\beta$, $C$ et $D$ soient 
$\Borel (\gamma)$, et posons 
$$F^{\gamma}_{\eta} := 
\left\{\!\!\!\!\!\!
\begin{array}{ll} 
& \overline{A' \cap \bigcap_{\theta<\eta} 
F^{\gamma}_{\theta}} ~\mbox{si}~\eta~\mbox{est~pair,} \cr & \cr
& \overline{\check A \cap \bigcap_{\theta<\eta} 
F^{\gamma}_{\theta}} ~\mbox{si}~\eta~\mbox{est~impair,}
\end{array}
\right.$$ 
l'adh\'erence \'etant prise au sens de la topologie 
${\it\Delta}^{\gamma}_X\times {\it\Delta}^{\gamma}_Y$. 
Alors il est clair que si $\eta \leq \xi$, 
$$F^{\gamma}_{\eta} \subseteq F_{\eta} \cap (C\times D)\mbox{,}$$ 
ce par r\'ecurrence transfinie. Il suffit donc de voir que 
$F^{\gamma}_{\xi}$ est non vide. En raisonnant par l'absurde, 
on va montrer que si $U_{\eta} = \check F^{\gamma}_{\eta}$, 
alors $A' = \bigcup_{\eta<\xi,\eta~\mbox{impair}} 
U_{\eta} \setminus (\bigcup_{\theta<\eta} U_{\theta})$. Si $\eta$ 
est impair et successeur de $\theta_0$, 
$U_{\eta} \setminus (\bigcup_{\theta<\eta} U_{\theta}) = 
F^{\gamma}_{\theta_0} \setminus F^{\gamma}_{\eta} = F^{\gamma}_{\theta_0} 
\setminus \overline{\check A \cap F^{\gamma}_{\theta_0}} \subseteq A'$. Si 
maintenant $x$ est dans $A'$, $x$ est dans $U_{\xi}$ donc il 
existe un plus petit $\eta \leq \xi$ tel que $x$ soit dans 
$U_{\eta}$. Si $\eta$ est pair et successeur de $\theta_0$, $x$ 
est dans $U_{\eta} \setminus (\bigcup_{\theta<\eta} U_{\theta}) 
= F^{\gamma}_{\theta_0} \setminus F^{\gamma}_{\eta} = F^{\gamma}_{\theta_0} 
\setminus \overline{A' \cap F^{\gamma}_{\theta_0}} \subseteq \check A'$. Si 
$\eta$ est limite, $x$ est dans 
$\check {\overline{A' \cap \bigcap_{\theta<\eta} F^{\gamma}_{\theta}}} \setminus 
(\bigcup_{\theta<\eta} U_{\theta}) = \check {\overline{A' \cap 
\bigcap_{\theta<\eta} F^{\gamma}_{\theta}}} \cap 
(\bigcap_{\theta<\eta} F^{\gamma}_{\theta}) \subseteq \check A'$. D'o\`u le 
r\'esultat.\bigskip

 On en d\'eduit que $F_{\xi}\cap (D ^{\beta}_X\times 
D ^{\beta}_Y)$ est un $\Ana (\beta)$ non vide, donc qu'il 
rencontre l'ensemble $\Omega ^{\beta}_{X\times Y} \subseteq \Omega 
^{\beta}_X\times \Omega ^{\beta}_Y$, et que 
$F_{\xi}\cap (Z_0\times T_0)$ est non vide. On remarque 
alors que dans la d\'efinition des $F_{\eta}$, pour 
$\eta \leq \xi$, on peut tout aussi bien prendre l'adh\'erence 
au sens de ${\it\Sigma}^{\beta}_X\times {\it\Sigma}^{\beta}_Y$ car 
les ensembles sous l'adh\'erence sont $\Ana (\beta)$. On a 
donc que $F_{\xi}\cap (Z_0\times T_0) = 
\overline{A \cap (Z_0\times T_0) \cap \bigcap_{\eta<\xi} 
F_{\eta}}^{Z_0\times T_0}$. L'ensemble $A \cap (Z_0\times T_0) \cap \bigcap_
{\eta<\xi} F_{\eta}$ est donc r\'eunion de $G_{\delta}$ 
l.p.o. de $Z_0\times T_0$, et le th\'eor\`eme 1.13 peut 
s'appliquer \`a l'un de ces $G_{\delta}$ 
qui est non vide ; ce qui fournit des presque-ouverts non 
vides $a_{\emptyset}$ et $b_{\emptyset}$, ainsi que 
$g_{\emptyset}:a_{\emptyset}\rightarrow b_{\emptyset}$ 
surjective continue ouverte.\bigskip

 On construit alors, pour $s$ dans $T_{\xi}$, des suites $(a_s)$ et $(b_s)$ de 
presque-ouverts (l'un dans $Z_0$, l'autre dans $T_0$), une 
suite $(g_s)$ de surjections continues ouvertes de $a_s$ sur 
$b_s$, des suites denses $(x^s_n)_n$ de $a_s$ v\'erifiant :
$$\begin{array}{ll}
& (i)\ ~G(g_s) \subseteq \left\{\!\!\!\!\!\!
\begin{array}{ll} 
& (Z_0\times T_0) \cap \bigcap_{n\in\omega} 
F_{f(s^{\frown} n)} \cap A ~\mbox{si}~ \vert s\vert ~
\mbox{est~paire}~((Z_0\times T_0)\cap A~{\rm si}~f(s) = 0)\mbox{,}\cr & \cr
& (Z_0\times T_0) \cap \bigcap_{n\in\omega} 
F_{f(s^{\frown} n)} \cap \check A ~\mbox{si}~ \vert s\vert ~\mbox{est~impaire,}
\end{array}
\right.\cr & \cr
& (ii)\ G(g_{s^\frown n})\mbox{ (ou }G(g_{s^\frown n})^*\mbox{) }\subseteq 
{\cal B}[(x^s_n,g_s(x^s_n)),2^{-\Sigma_{i \leq \vert s\vert }~(s^\frown n(i)+1)}]\mbox{ si }s \in T'_{\xi}.
\end{array}$$
Admettons avoir construit $g_s$ pour $\vert s\vert  \leq p$, 
et soient $s$ de longueur $p$, et $(x^s_n)_n$ une suite 
dense de $a_s$. Alors si $p$ est pair et $n$ entier, $G(g_s) 
\subseteq F_{f(s^{\frown} n)}\cap (Z_0\times T_0)$, qui est \'egal \`a 
$$(Z_0\times T_0) \cap \overline{\bigcap_{m\in\omega} 
F_{f(s^{\frown} n^{\frown} m)}\setminus A}.$$

\vfill\eject

 D'o\`u 
$G(g_s) \subseteq \overline{(Z_0\times T_0)\cap 
\bigcap_{m\in\omega} F_{f(s^{\frown} n^{\frown} m)}
\setminus A}$. 
On peut alors appliquer le th\'eor\`eme 1.13 \`a l'un des ensembles
$G_{\delta}$ et $\Ana (\beta)$ non vides dont l'intersection 
$$\bigcap_{m\in\omega} F_{f(s^{\frown} n^{\frown} m)}
\setminus A \cap {\cal B}
[(x^s_n,g_s(x^s_n)), 2^{-\Sigma_{i \leq \vert s\vert }~(s^\frown n(i)+1)}]$$ 
(ou alors l'intersection $\bigcap_{m\in\omega} F_{f(s^{\frown} n^{\frown} m)}
\setminus A\cap {\cal B}
[(g_s(x^s_n),x^s_n), 2^{-\Sigma_{i \leq \vert s\vert }~(s^\frown n(i)+1)}]$) 
est la r\'eunion, dans le produit $Z_0\times T_0$, ce qui fournit 
le graphe recherch\'e. De m\^eme si $p$ est impair.\bigskip
 
 On choisit alors $\alpha$ tel que $Z_0$ et $T_0$ soient 
r\'ecursivement en $\alpha$-pr\'esent\'es, et tel que pour tout $s\in T_\xi$ et pour 
tout $p\in \omega$, $\bigcup_{n\in\omega} G(g_{s^\frown n})$, $G(g_s)$ et 
$\bigcup_{s\in \omega^{\leq p}} G(g_s)$ soient $\Borel (\alpha)$. On prend des 
notations analogues aux pr\'ec\'edentes. Soit $M_0$ (resp. $N_0$) un $G_{\delta}$ 
dense $\Ana(\alpha)$ de $Z_0$ (resp. $T_0$), fourni par le lemme 2.1, sur 
lequel la topologie de $Z_0$ (resp. $T_0$) co\"\i ncide avec la topologie 
${\it\Delta}^{\alpha}_{Z_0}$ (resp. ${\it\Delta}^{\alpha}_{T_0})$. On 
d\'efinit maintenant les objets recherch\'es :
$$Z~:=~M_0~\cap~\Omega^{\alpha}_{Z_0}~\cap ~D^{\alpha}_{Z_0}\mbox{, }
T~:=~N_0~\cap~\Omega^{\alpha}_{T_0}~\cap ~D^{\alpha}_{T_0}\mbox{,}$$ 
$$G(f_s) := \left \{\!\!\!\!\!\!
\begin {array}{ll} 
& G(g_s) \cap 
(Z\times T)~\mbox{si}~f(s) = 0\mbox{,}\cr & \cr
& G(g_s) \cap \overline{C_s}^{Z\times T}
~\mbox{si}~s \in T'_{\xi}.
\end{array}
\right.$$
$A_s$ et $B_s$ sont les projections de $G(f_s)$, et $u$ (resp. 
$v$) est l'application identique de $Z$ dans $X$ (resp. $T$ dans $Y$). 
V\'erifions que ces objets conviennent. Montrons que $G(f_s) = G(g_s) \cap (Z\times T)$ si $s \in T_{\xi}$. \bigskip

 La relation est vraie par d\'efinition si $f(s) = 0$. 
Admettons-la pour $f(s)<\eta \leq \xi$ ; soit $s \in 
T_{\xi}$ tel que $f(s) = \eta > 0$. On a bien s\^ur que 
$G(f_s) \subseteq G(g_s)\cap (Z\times T)$. Mais on a 
$$G(f_s) = G(g_s)\cap \overline{\bigcup_{n\in\omega} 
G(f_{s^\frown n}})^{Z\times T} = {G(g_s)\cap \overline{
\bigcup_{n\in\omega} G(g_{s^\frown n})\cap (Z\times T)}^{Z\times T}}\mbox{,}$$ 
par hypoth\`ese de r\'ecurrence. Montrons donc que 
$G(g_s)\cap (Z\times T)\subseteq\overline{\bigcup_{n\in\omega} G(g_{s^\frown n})\cap 
(Z\times T)}^{Z\times T}$. On a $G(g_s) \subseteq \overline{
\bigcup_{n\in\omega} G(g_{s^\frown n})}^{Z_0\times T_0}$ car si 
$O$ est ouvert-ferm\'e et contient $(x,g_s(x))$, on trouve une 
suite strictement croissante $(n_q)_q$ telle que 
$(x^s_{n_q})_q$ converge vers $x$. Par continuit\'e de $g_s$, 
la suite image converge vers $g_s(x)$, et on trouve $q$ tel 
que ${{\cal B} [(x^s_{n_q},g_s(x^s_{n_q})), 2^{-\Sigma_{i \leq 
\vert s\vert }~(s^\frown n_q(i)+1)}] \subseteq O}$, et $G(g_{s^\frown n_q}) 
\subseteq O$ (ou $G(g_{s^\frown n_q}) \subseteq O^*$). On a que 
$\Omega^{\alpha}_{Z_0}$ est co-maigre dans $Z_0$. En effet, $Z_0$ est 
polonais parfait de dimension 0 et non vide, donc on peut choisir 
l'isomorphisme avec $\omega^\omega$ de fa\c con \`a pr\'eserver les ensembles co-maigres ; 
il suffit alors de consulter [Ke], o\`u il est d\'emontr\'e que 
$\Omega^{\alpha}_{\omega^\omega}$ est co-maigre dans $\omega^\omega$. Donc $Z$ 
(resp. $T$) est co-maigre dans $Z_0$ (resp. $T_0$) et 
$$G(g_s)\cap (Z\times T) \subseteq \overline{
\bigcup_{n\in\omega} G(g_{s^\frown n})\cap (Z\times T)}^
{Z_0\times T_0}\cap (Z\times T).$$ 
Or ce dernier ensemble est $\overline{ \bigcup_{n\in\omega} G(g_{s^\frown n})\cap 
(Z\times T)}^{Z\times T}$, car sur $Z$ (resp. $T$), les 
topologies initiales et ${\it\Delta}^{\alpha}_{Z_0}$ (resp. 
${\it\Delta}^{\alpha}_{T_0}$) co\"\i ncident, et car on a un $\Ana (\alpha )$ 
sous l'adh\'erence, donc les adh\'erences pour 
${\it\Delta}^{\alpha}_{Z_0} \times {\it\Delta}^{\alpha}_{T_0}$ et 
${\it\Sigma}^{\alpha}_{Z_0} \times {\it\Sigma}^{\alpha}_{T_0}$ sont les 
m\^emes.

\vfill\eject

 Comme $\Omega^{\alpha}_{Z_0}\cap D^{\alpha}_{Z_0}$, muni de 
la restriction de la topologie de Gandy-Harrington, est 
polonais parfait de dimension 0, $Z$ l'est aussi puisqu'il 
en est un ouvert. De m\^eme pour $T$. Les seules choses non 
\'evidentes \`a v\'erifier sont que $G(f_s)$ n'est pas vide, l'inclusion    
${\overline{C_s} \setminus C_s \subseteq G(f_s)}$, et aussi que 
$\overline{B_p} = B_p \cup B_i$.\bigskip

 L'ensemble $G(f_s)$ est non vide puisqu'il est \'egal \`a 
$G(g_s)\cap (Z\times T)$ et que $Z$ (resp. $T$) est co-
maigre dans $Z_0$ (resp. $T_0$).\bigskip

  Soit $(x,y)$ dans $\overline{C_s}^{Z\times T} 
\setminus C_s$, $((x_n,y_n))$ dans $C_s$ convergeant vers 
$(x,y)$, et $p_n$ tel que $(x_n,y_n)$ soit dans $G(f_{s^\frown p_n})$. 
Comme $G(f_{s^\frown p_n})$ est ferm\'e dans $Z\times T$, on 
peut supposer la suite $(p_n)$ strictement croissante. La 
distance de $(x_n,y_n)$ \`a $G(g_s)$, dans $Z_0\times T_0$, 
est au plus $2^{-\Sigma_{i\leq \vert s\vert }~(s^\frown p_n(i)+1)}$, 
et comme la convergence a lieu aussi dans $Z_0\times T_0$, 
on a $(x,y) \in \overline{G(g_s)}^{Z_0\times T_0}$. \bigskip

 Mais on a l'\'egalit\'e
$\overline{G(g_s)}^{Z_0\times T_0} \cap (Z\times T) = 
\overline{G(g_s) \cap (Z\times T)}^{Z_0\times T_0} \cap 
(Z\times T)$ car $Z$ et $T$ sont co-maigres, d'o\`u $(x,y) 
\in \overline{G(f_s)}^{Z\times T}$, comme pr\'ec\'edemment. Comme 
le graphe de $f_s$ est ferm\'e, on a bien $(x,y) \in G(f_s)$.\bigskip

 Montrons par r\'ecurrence sur $p$ que $G(f_{\emptyset}) 
\cup \bigcup_{\vert s\vert <p,s\in T'_{\xi}} C_s$ est ferm\'e. 
C'est clair pour $p=0$. Admettons donc que c'est vrai pour 
 $p$. On a les \'egalit\'es 
 $$G(f_{\emptyset})\cup 
\!\!\bigcup\limits_{\vert s\vert <p+1}\!\!C_s =
G(f_{\emptyset})~\cup\!\!\bigcup\limits_{\vert s\vert <p}\!\! C_s \cup 
\!\!\bigcup\limits_{s\in \omega^p}\!\! C_s=
{G(f_{\emptyset}) \cup\!\! 
\bigcup\limits_{\vert s\vert <p}\!\! C_s \cup \!\!\bigcup\limits_{s\in 
\omega^{p +1}}\!\! G(f_s)}.$$ 
Soit $((x_m,y_m)) \subseteq G(f_{\emptyset}) \cup 
\bigcup_{\vert s\vert <p+1} C_s$ convergeant vers $(x,y)$. Alors 
par hypoth\`ese de r\'ecurren-ce, on peut supposer que pour chaque 
$m$ il existe $(s_m,n_m)$ tel que $(x_m,y_m)$ soit dans 
$G(f_{s_m^\frown n_m})$. Les ensembles $G(f_{s^\frown n})$ \'etant ferm\'es, on peut 
supposer que la suite $(\Sigma_{i\leq p}~(s_m^\frown n_m(i)+1))_m$ 
tend vers l'infini. La distance de $(x_m,y_m)$ \`a $G(g_{s_m})$, 
dans $Z_0\times T_0$, est au plus $2^{-\Sigma_{i\leq p}~(s_m^\frown n_m(i)+1)}$, 
donc comme la convergence a lieu aussi dans $Z_0\times T_0$, on a 
$$(x,y) \in \overline{G(g_{\emptyset}) 
\cup \bigcup_{(s,n)\in \omega^{<p}\times \omega} 
G(g_{s^\frown n})}^{Z_0\times T_0} \cap (Z\times T).$$ 
Mais en raisonnant comme pr\'ec\'edemment, on voit que 
$$(x,y) \in\overline{G(f_{\emptyset}) \cup \bigcup_{s\in \omega^{<p}}
 C_s}^{Z\times T} = G(f_{\emptyset}) \cup \bigcup_{s\in 
\omega^{<p}} C_s\mbox{,}$$ 
par hypoth\`ese de r\'ecurrence. Donc $(x,y) \in G(f_{\emptyset}) \cup \bigcup_{s\in \omega^{<p+1}} 
C_s$, qui est donc ferm\'e.\bigskip

 On a que $G(f_{\emptyset}) \cup \bigcup_{s\in T'_{\xi}} C_s$ est ferm\'e. En effet, 
si $((z_m,t_m)) \subseteq G(f_{\emptyset}) \cup \bigcup_{s\in T'_{\xi}} C_s$ converge 
vers $(z,t)$, par ce qui pr\'ec\`ede on peut supposer que pour chaque $m$ il existe 
$s'_m$ tel que $(z_m,t_m)$ soit dans $G(f_{s'_m})$, et que la 
suite $(\vert s'_m\vert )_m$ tend vers l'infini. Alors on trouve $p$ tel que l'ensemble 
des $s'_m(p)$ soit infini.

\vfill\eject

 En effet, si tel n'est pas le cas, $\{ s\in T_\xi~/~
\exists~m~~s\prec s'_m\}$ est un sous-arbre infini de $T_\xi$, \`a branchements finis, 
donc a une branche par le lemme de K\"onig ; mais ceci contredit la bonne fondation 
de $T_\xi$. On peut donc supposer qu'il existe $p$ tel que la suite 
$(s'_m\lceil p)_m$ soit constante, disons \`a $s$, et tel que la suite $(s'_m(p))_m$ tende vers 
l'infini. Alors par l'in\'egalit\'e triangulaire on a que la distance de $(z_m,t_m)$ 
\`a $G(g_s)$, dans $Z_0\times T_0$, est au plus $2^{-s'_m(p)}$, donc comme la convergence 
a lieu aussi dans $Z_0\times T_0$, on a $(z,t) \in \overline{G(g_s)}^
{Z_0\times T_0}\cap (Z\times T)$. Mais comme avant, on en d\'eduit que 
$(z,t)$ est dans $G(f_s)\subseteq G(f_{\emptyset}) \cup \bigcup_{s\in T'_{\xi}} C_s$.\bigskip

 On a donc que $\overline{B_p} \subseteq G(f_{\emptyset}) \cup 
\bigcup_{s\in T'_{\xi}} C_s$. Si $\vert s\vert $ est paire, 
on montre que $C_s \subseteq \overline{B_p}$, ce qui 
montrera que $\overline{B_p} = G(f_{\emptyset}) \cup \bigcup_
{s\in T'_{\xi}} C_s = B_p \cup B_i$. Or on a 
$$C_s = \bigcup_{n\in\omega} 
G(f_{s^\frown n}) \subseteq \bigcup_{n\in\omega} 
\overline{C_{s^\frown n}} \subseteq \overline{\bigcup_{n\in
\omega} C_{s^\frown n}} \subseteq \overline{B_p}.$$
Ceci termine la preuve.$\hfill\square$\bigskip
 
 Par [Lo-SR], on a, pour $\xi<\omega_1$ impair, l'existence d'un 
compact $K_\xi$ et d'un vrai $D_\xi (\boraone)$ de $K_\xi$, disons $B_\xi$, 
tel que si $A$ est bor\'elien d'un espace polonais $P$, $A$ n'est pas 
$\check D_\xi (\boraone)$ de $P$ si et seulement s'il existe $f:K_\xi 
\rightarrow P$ injective continue telle que $f^{-1}(A) = B_\xi$. On a 
ici un analogue : si $A$ est bor\'elien d'un produit d'espaces polonais, $A$ 
est non-$\mbox{pot}(\check D_\xi (\boraone))$ si et seulement s'il existe 
$B$ dans $D_\xi (\boraone)\lceil \overline{B} \setminus \mbox{pot}(\check D_\xi (\boraone))$ 
et des injections continues $u$ et $v$ tels que $B = \overline{B}\cap 
(u\times v)^{-1}(A)$. Bien s\^ur, $B$ d\'epend ici de $A$, mais est toujours du 
m\^eme type $(B = D((\bigcup_{f(s)\leq \eta} G(f_s))_{\eta<\xi}))$.\bigskip

 On va maintenant donner des versions du th\'eor\`eme 2.3 
dans le cas o\`u $A$ est \`a coupes verticales d\'enombrables, et 
dans le cas o\`u $\xi=1$.

\begin{lem} (a) Soient $X$ et $Y$ des espaces polonais, 
$A$ (resp. $B$) un bor\'elien de $X$ (resp. $Y$), et 
$f:A\rightarrow B$ countable-to-one bor\'elienne. Alors il 
existe une partition de $A$ en ensembles bor\'eliens, disons
 $(A_n)_{n\in \omega}$, telle que les restrictions de $f$ \`a $A_n$ soient 
injectives.\smallskip

\noindent (b) Soient $X$ et $Y$ des espaces polonais,
 $A$ (resp. $B$) un presque-ouvert non vide de $X$ (resp. $Y$), 
et $g:B\rightarrow A$ surjective continue ouverte 
countable-to-one. Alors il existe un presque-ouvert non vide 
 $A'$ (resp. $B'$) de $X$ (resp. $Y$), contenu dans $A$ (resp. 
$B$), tel que la restriction de $g$ \`a $B'$ soit un hom\'eomorphisme 
de $B'$ sur $A'$.\end{lem}

\noindent\bf D\'emonstration.\rm\ (a) Par le th\'eor\`eme de Lusin, 
l'ensemble $\mbox{Gr}(f)^* := \{ (y,x)~/~(x,y) \in \mbox{Gr}(f) \}$ est la r\'eunion
 d\'enombrable des graphes de fonctions bor\'eliennes partielles, 
disons $f_n$, d\'efinies sur des bor\'eliens. Puisque leur 
graphe est contenu dans $\mbox{Gr}(f)^*$, les $f_n$ sont injectives, 
donc leurs images $C_n$ sont bor\'eliennes. Il reste \`a poser 
 $A_n := C_n \setminus \bigcup_{p<n} C_p$.\bigskip

\noindent (b) Soit $(B_n)$ une partition bor\'elienne de $B$ 
donn\'ee par le (a). Comme $B$ est non maigre relativement
 \`a $Y$, l'un des $B_n$ est non maigre. Ce
 dernier ayant la propri\'et\'e de Baire s'\'ecrit comme r\'eunion 
disjointe d'un $G_{\delta}$ non maigre $G$ de $Y$ et d'une
 partie maigre $M$. Soit $U$ un ouvert non vide de $Y$ tel que 
 $G \Delta U$ soit maigre, et posons $B'' := U\cap G$. Cet
 ensemble est un $G_{\delta}$ dense de $B\cap U$, qui est 
ouvert de $B$, donc $g[B'']$ est un analytique co-maigre de 
 l'ouvert $g[B\cap U]$ de $A$.
 
\vfill\eject
 
  Soient donc $(U_n)$ une 
base de la topologie de $B''$, $V_n$ des ouverts de 
$g[B\cap U]$ tels que $g[U_n] \Delta V_n$ soit un maigre $M_n$ 
dans $g[B\cap U]$, et $A'$ un $G_{\delta}$ dense de $g[B\cap U]$ 
contenu dans $g[B'']\setminus \bigcup_{n\in\omega} M_n$. 
Alors on peut poser $B' := B'' \cap g^{-1}(A')$, comme on le 
v\'erifie facilement.$\hfill\square$

\begin{thm} Soient $X$ et $Y$ des espaces 
polonais, $A$ un bor\'elien \`a coupes verticales d\'enombrables de 
$X \times Y$, et $\xi$ un ordinal d\'enombrable non nul.\smallskip

\noindent (a) Si $\xi$ est pair, $A$ est non-$\mbox{pot}(D_{\xi}(\boraone))$ si et seulement 
s'il existe des espaces polonais parfaits $Z$ et $T$ de dimension 0, des 
ouverts-ferm\'es non vides $A_s$ et $B_s$ (l'un dans $Z$ et l'autre dans $T$, 
pour $s$ dans $T_{\xi}$, et dans cet ordre si $\vert s\vert $ est paire), 
des surjections continues ouvertes $f_s$ de $A_s$ sur $B_s$, 
et des injections continues $u$ et $v$ tels que si 
$C_s := \bigcup_{n\in\omega} G(f_{s^\frown n})$, on ait 
$$G(f_s) \subseteq \overline{C_s}\setminus (\bigcup_{t \in T'_{\xi}~/~
\mbox{parit\'e}(\vert t\vert )=\mbox{parit\'e}(\vert s\vert )} C_t)\mbox{,}$$ 
et si $B := \bigcup_{s\in T_{\xi}~/~\vert s\vert ~\mbox{paire}} G(f_s)$, 
$B = (u\times v)^{-1}(A)$.\smallskip

\noindent (b) Si $\xi$ est impair, $A$ est non-$\mbox{pot}(\check D_{\xi}(\boraone))$ si et 
seulement s'il existe des espaces polonais $Z$ et $T$ parfaits de 
dimension 0, des ouverts-ferm\'es non vides $A_s$ 
et $B_s$ (l'un dans $Z$ et l'autre dans $T$, pour $s$ dans 
$T_{\xi}$, et dans cet ordre si $\vert s\vert $ est 
impaire), des surjections continues ouvertes 
$f_s$ de $A_s$ sur $B_s$, et des injections continues $u$ 
et $v$ tels que si 
$C_s := \bigcup_{n\in\omega} G(f_{s^\frown n})$, on ait 
$${G(f_s) \subseteq \overline{C_s}\setminus (\bigcup_{t \in T'_{\xi}~/~
\mbox{parit\'e}(\vert t\vert )=\mbox{parit\'e}(\vert s\vert )} C_t)}\mbox{,}$$ 
et si $B := \bigcup_{s\in T_{\xi}~/~\vert s\vert ~\mbox{impaire}} G(f_s)$, on ait
${B = (u\times v)^{-1}(A)}$.\end{thm}

\noindent\bf D\'emonstration.\rm\ Elle est identique \`a la preuve du th\'eor\`eme 
 2.3, \`a ceci pr\`es qu'on remplace la condition (ii) par 
 $\mbox{Gr}(g_{s^\frown n})$ (ou $\mbox{Gr}(g_{s^\frown n})^*$) $\subseteq 
 {\cal B}[(x^s_n,g_s(x^s_n)),2^{-\Sigma_{i \leq \vert s\vert } 
 ~(s^\frown n(i)+1)}]$ si $s\in T'_{\xi}$ et $f(s) \not= 1$, et $A\cap 
 (Z_0\times T_0) = \bigcup_{n\in\omega} \mbox{Gr}(g_{s^\frown n})$ si 
$f(s) = 1$. Si $A$ est la $\Borel (\beta)$-r\'eunion de $(h_n)$, 
on prend pour la suite $(\mbox{Gr}(g_{s^\frown n}))_n$ la suite des 
traces de $\mbox{Gr}(h_n)$ sur $Z_0\times T_0$ qui sont non vides 
(ceci si $f(s) = 1$). Si $f(s)>0$ est pair, 
on raisonne comme dans la preuve de 2.3 pour construire $g_s$, 
\`a ceci pr\`es que si le graphe n'est pas contenu dans $Z_0\times T_0$ 
mais dans $T_0\times Z_0$, on applique le lemme 2.4 pour avoir 
le graphe dans le sens souhait\'e.$\hfill\square$

\begin{thm} Soient $X$ et $Y$ des espaces 
polonais, $A$ un bor\'elien \`a coupes verticales d\'enombrables de 
$X \times Y$. Alors $A$ est non-$\mbox{pot}(\bormone)$ si et seulement s'il 
existe des espaces polonais $Z$ et $T$ parfaits de dimension 0, 
une suite d'ouverts-ferm\'es ($A_n$) (resp. ($B_n$)) de $Z$ (resp. $T$), 
des surjections continues ouvertes $f_n$ de $A_n$ sur $B_n$, et des 
fonctions continues $u$ et $v$ tels que si $B = \bigcup_{n>0} \mbox{Gr}(f_n)$, 
on ait $\emptyset \not= \mbox{Gr}(f_0) \subseteq \overline{B} \setminus B$ et 
$B=(u\times v)^{-1}(A)$.\end{thm}

\noindent\bf D\'emonstration.\rm\ Elle est identique \`a celle du th\'eor\`eme 2.5 si 
$a_{\emptyset} \subseteq Z_0$. Dans le cas o\`u $a_{\emptyset} \subseteq T_0$, 
on montre comme dans 2.5 que l'ensemble 
$$L:=G(g_{\emptyset}) \cap  [(M_0 \cap \Omega_{Z_0}^{\alpha} \cap D_{Z_0}^{\alpha})
\times  (N_0 \cap \Omega_{T_0}^{\alpha} \cap D_{T_0}^{\alpha})] \cap 
\overline{A\cap (Z_0\times T_0)}^{{\it\Sigma}_{Z_0}^{\alpha}\times {\it\Sigma}_{T_0}^{\alpha}}$$ 
est non vide.

\vfill\eject

 On d\'efinit alors $A_0$ 
comme \'etant la projection de $L$ sur $T_0$, $B_0 := Z := T := A_0$ ; 
$f_0$ est l'application identique. Soient $(h_n)_{n>0}$ comme dans la 
preuve de 2.5, $E_n$ le domaine de $h_n$, $D_0$ la projection de $L$ sur 
$Z_0$, et $g_0:A_0\rightarrow D_0$ la restriction de $g_{\emptyset}$ \`a
 $A_0$. On d\'efinit $A_n := g_0^{-1}(h_n^{-1}(A_0))$, 
 $B_n := A_0\cap h_n[E_n\cap D_0]$, et $f_n$ comme \'etant la restriction de 
 $h_n\circ g_0$ \`a $A_n$, si $n>0$. Enfin, on pose $u:=g_0$, $v(y)=y$ si 
 $y\in A_0$. V\'erifions que ces objets conviennent. 
Si $x\in A_0$, $(x,x) \notin B$ sinon $(g_0(x),x) \in A\cap 
  \mbox{Gr}(g_{\emptyset})^*$. Si $U$ est un ouvert de $A_0$ contenant $x$, 
  $g_0[U]\times U$ est un ouvert de $(M_0 \cap \Omega_{Z_0}^
  {\alpha} \cap D_{Z_0}^{\alpha})\times  (N_0 \cap \Omega_{T_0}^{\alpha} 
  \cap D_{T_0}^{\alpha})$ contenant $(g_0(x),x)$, donc rencontre $A$ en un point 
  $(z,y)$ ; $z = g_{\emptyset}(t)$, o\`u $t\in U$, et $(t,y)\in U^2\cap B \not= 
  \emptyset$. Donc ${\Delta (A_0) = \mbox{Gr}(f_0) \subseteq \overline{B} \setminus B}$.\bigskip 
  
Enfin, $(x,y) \in B \Leftrightarrow \exists~n~~y = h_n(g_0(x)) \Leftrightarrow 
\exists~n~~(g_0(x),y) \in \mbox{Gr}(h_n) 
\Leftrightarrow (g_0(x),y) \in A$.$\hfill\square$\bigskip 

 Soit $C_0$ la classe des fonctions de la forme  
$(x,y)\!\mapsto \! (u(x),v(y))$ ou 
${(x,y)\!\mapsto \! (u(y),v(x))}$, o\`u les 
fonctions $u$ et $v$ sont bor\'eliennes.

\begin{thm} Soient $X$ et $Y$ des espaces 
polonais, $A$ un bor\'elien $\mbox{pot}(\borathree)\cap \mbox{pot}(\bormthree)$ de 
$X \times Y$. \smallskip

\noindent (a) $A$ est non-$\mbox{pot}(\bormone)$ si et seulement s'il 
existe des espaces polonais $Z'$ et $T'$ parfaits de dimension 0, 
une suite d'ouverts-ferm\'es ($A_n$) (resp. ($B_n$)) de $Z'$ (resp. $T'$), 
des surjections continues ouvertes $f_n$ de $A_n$ sur $B_n$, et une 
fonction continue $f$ de $C_0$ tels que si ${B = \bigcup_{n>0} \mbox{Gr}(f_n)}$, 
on ait 
$$\emptyset \not= \mbox{Gr}(f_0) = \overline{B} \setminus B$$ 
et $B=f^{-1}(A)\cap \overline{B}$.\medskip

\noindent (b) $A$ est non-$\mbox{pot}(\bormone)$ si et seulement s'il 
existe des espaces polonais $Z$ et $T$ parfaits de dimension 0 non vides, 
une suite d'ouverts denses $(E_n)$ de $Z$, 
des applications continues ouvertes $g_n$ de $E_n$ dans $T$, et une 
fonction continue $g$ de $C_0$ tels que $(g_n)_{n>0}$ converge uniform\'ement vers 
$g_0$ sur $\bigcap_{n\in\omega} E_n$, $\mbox{Gr}(g_0) \subseteq g^{-1}(\check A)$ et 
$\bigcup_{n>0} \mbox{Gr}(g_n) \subseteq g^{-1}(A)$.\end{thm}

\noindent\bf D\'emonstration.\rm\ (a) Dans le cas o\`u $a_{\emptyset} \subseteq Z_0$, la preuve 
 est identique \`a celle du th\'eor\`eme 2.3 jusqu'\`a la construction des $g_s$ 
 comprise. Soient $I_0 := \{ n~/~a_{(n)} \subseteq Z_0 \}$, 
 ${I_1 := \{ n~/~a_{(n)} \subseteq T_0 \}}$, 
 $$J_0 := \overline{\bigcup_{n\in I_0} \mbox{Gr}(g_{(n)})}^{Z_0\times T_0}\mbox{,}$$ 
 $J_1 := \overline{\bigcup_{n\in I_1} \mbox{Gr}(g_{(n)})^*}^{Z_0\times T_0}$. Alors 
 $\mbox{Gr}(g_{\emptyset}) \subseteq J_0 \cup J_1$, donc on trouve $i$ et un 
 ouvert-ferm\'e $U\times V$ de $Z_0\times T_0$ tels que 
 $\emptyset \not= \mbox{Gr}(g_{\emptyset}) \cap (U\times V) \subseteq J_i$. Supposons 
 par exemple que $i=0$. Soit $g_0$ la restriction de $g_{\emptyset}$ \`a 
 $U\cap g_0^{-1}(V)$, et $B' := \bigcup_{n\in I_0} \mbox{Gr}(g_{(n)})$ ; on a donc 
 $\mbox{Gr}(g_0) \subseteq \overline{B'}^{Z_0\times T_0}$. On conclut de mani\`ere 
 analogue \`a celle de 2.3 : on introduit $\alpha$, on montre que l'ensemble
 $$K := \mbox{Gr}(g_0) \cap (Z\times T) \cap \overline{B'\cap [(Z\cap U)\times 
 (T\cap V)]}^{{\it\Sigma}_{Z_0}^{\alpha}\times {\it\Sigma}_{T_0}^{\alpha}}$$ 
 est non vide, ce qui implique que $I_0$ est infini ; on indexe alors $B'$ par $\omega$. On d\'efinit 
 $A_0$ et $B_0$ (resp. $A_n$ et $B_n$) comme \'etant les projections de $K$ 
 (resp. $\mbox{Gr}(g_{(n)})\cap [(Z\cap U)\times (T\cap V)]$), et on pose 
 $f_0 := g_0 \lceil A_0$, ${f_n := g_{(n)} \lceil A_n}$ si $n>0$, 
$Z' := Z\cap U$, $T' := T\cap V$, et $f$ est l'application identique.

\vfill\eject

  Si maintenant $i=1$, on montre comme dans 2.6 qu'on peut trouver deux 
  fonctions continues $u$ et $v$, un espace polonais $Z'$ parfait de 
  dimension 0, et des applications continues et ouvertes d\'efinies sur et 
  d'images des ouverts-ferm\'es de $Z'$, $\varphi_n$, tels que si 
  $A' := \bigcup_{n>0} \mbox{Gr}(\varphi_n)$, on ait 
  $$\Delta (Z') \subseteq (v\times u)^{-1}(\mbox{Gr}(\varphi_0))\mbox{,}$$ 
  $\mbox{Gr}(\varphi_0) \subseteq 
  \mbox{Gr}(g_{\emptyset})\cap (U\times V)$, $A'^* \subseteq (v\times u)^{-1}(A)$, et 
  $\Delta (Z') \subseteq \overline{A'} \setminus A'$. Soit $(t_n)$ une suite 
  dense de $Z'$. Alors on construit des rectangles ouverts-ferm\'es $W_n$ de 
  $Z'^2$ de diam\`etre au plus $2^{-n}$ et une suite strictement croissante 
  $(m_n)$ telle que l'on ait $(t_n,t_n) \in W_n$, 
  $$\emptyset \not= \mbox{Gr}(\varphi_{m_n})\cap W_n \subseteq ~^c(\bigcup_{p<n} 
\mbox{Gr}(\varphi_{m_p})).$$
Il est maintenant clair qu'on peut poser $A_0 := B_0 := T' := Z'$, $f_0 := \mbox{Id}_{A_0}$ ; si $n>0$, $A_n$ et $B_n$ sont les projections de $\mbox{Gr}(\varphi_{m_n})\cap W_n$, et $f_n$ est la restriction de 
$\varphi_{m_n}$ \`a $A_n$ ; enfin, 
$$f(x,y) := (v(y),u(x)).$$
Le cas o\`u $a_{\emptyset} \subseteq T_0$ se traite de fa\c con analogue.\bigskip

\noindent (b) Supposons que $A$ est $\mbox{pot}(\bormone)$, en raisonnant par l'absurde. Alors 
si on pose ${G := \bigcap_{n\in\omega} E_n}$, $\bigcup_{n\in\omega} \mbox{Gr}(g_n) 
\cap (G\times T)$ est bor\'elien \`a coupes ferm\'ees, donc on trouve un $G_\delta$ 
dense $H\subseteq G$ tel que $\bigcup_{n\in\omega} \mbox{Gr}(g_n) \cap (H\times T)$ 
soit ferm\'e dans $H\times T$. Par ailleurs, 
$$g^{-1}(A)\cap \bigcup_{n\in\omega} \mbox{Gr}(g_n) \cap (H\times T) = 
\bigcup_{n>0} \mbox{Gr}(g_n) \cap (H\times T)$$ 
est $\mbox{pot}(\bormone)$, donc on trouve un $G_\delta$ dense $K\subseteq T$ tel que les coupes de
$E := \bigcup_{n>0} \mbox{Gr}(g_n) \cap (H\times K)$ soient ferm\'ees dans $K$. Alors 
$H \cap \bigcap_{n\in\omega} g_n^{-1}(K)$ est un $G_{\delta}$ 
dense de $Z$, et si $x$ est dans ce $G_{\delta}$, $g_0(x)$ est dans 
$\overline{E_x} \setminus E_x$, au sens de $K$, ce qui est la contradiction cherch\'ee.\bigskip

 Inversement, soient $Z'$, $T'$, $A_n$, $B_n$, $f_n$, $f$, fournis par le (a). 
Posons, si $n>0$, 
$${H_{\!n} \!\!:= \!\!\{ P \!\!\subseteq\! \!\borone \!\lceil\! A_0
\! \!\setminus \!\!\{ \!
\emptyset \!\}\!~/~\!\forall (\!U,V\!)\! \!\in\! \!P^2~(\!U\!\!\not=\!\!V\! 
\!\Rightarrow \! U 
\!\cap \!V\!\!\! =\!
\!\emptyset)~ \mbox{et} \exists p\!\!>\!\!0~(\!U \!\!\subseteq \!\!A_{\!p} ~\mbox{et}~ 
\forall x \!\!\in\!\! U~d(\!f_{\!0}(\!x\!),f_{\!p}(\!x\!)\!)\!<\! \!2^{-n}\!) \}}.$$
Alors $H_n$ est non vide, puisqu'il contient $\emptyset$, et il est 
ordonn\'e de mani\`ere inductive par l'inclusion, donc il admet un \'el\'ement maximal 
$P_n := \{ U^n_m~/~m \in I_n \}$. Posons ${E_n := \bigcup_{m\in I_n} U^n_m}$ ; 
$E_n$ est dense dans $A_0$, par maximalit\'e de $P_n$. On pose alors $Z := A_0$, 
$T := T'$, $g := f\lceil (Z\times T)$, $E_0 := A_0$, $g_0 := f_0$, et si $n>0$, 
$g_n (x) := f_{p_m} (x)$ si $x \in U^n_m$, o\`u bien s\^ur $p_m$ est minimal tel que 
$U^n_m \subseteq A_{p_m}$ et $\forall~x \in U^n_m~~d(f_0(x), f_{p_m}(x)) < 
2^{-n}$. Ces objets r\'epondent clairement au probl\`eme.$\hfill\square$
   
\begin{lem} Soient $X$ et $Y$ des espaces polonais, $A$ un bor\'elien 
  \`a coupes verticales d\'enombrables de $X \times Y$, et $k$ un entier naturel 
  non nul. Alors $\Gamma_A \not= D_k (\boraone)^+$, $\check D_{2k} (\boraone)$, et 
  $D_{2k-1} (\boraone)$ ($\Gamma_A$ est la plus petite classe de Wadge $\Gamma$ 
telle que $A\in \mbox{pot}(\Gamma)$, pour l'inclusion - cf [Le1] ; 
$\Gamma^+$ est le successeur de $\Gamma$).\end{lem} 

\vfill\eject
  
\noindent\bf D\'emonstration.\rm\ Elle repose sur le fait que si $B$ est bor\'elien \`a coupes 
   verticales d\'enombrables et est $\mbox{pot}(\boraone)$, $B$ est en fait $\mbox{pot}(\borone)$ 
   (puisque sa projection est d\'enombrable).\bigskip
   
    Il suffit de montrer que $\Gamma_A \not= \check D_{2k} (\boraone)$, la preuve 
    de la troisi\`eme assertion \'etant analogue et celle de la premi\`ere se d\'eduisant 
    des deux autres.\bigskip
    
 Supposons que 
$A = \check U_{2k-1} \cup U_0 \cup (U_2\setminus U_1) \cup ... \cup 
 (U_{2k-2}\setminus U_{2k-3})$ ; on a alors aussi 
$$A = ~^c(U_{2k-1}\setminus U_0) \cup 
 [(U_2\setminus U_0)\setminus (U_1\setminus U_0)] \cup ... \cup 
[(U_{2k-2}\setminus U_0)\setminus (U_{2k-3}\setminus U_0)].$$ 
Donc si on pose 
 $V_{2k-1} := X\times Y$ et $V_p := U_{p+1}\setminus U_0$ si $p<2k-1$, on a 
$$A = (V_1\setminus V_0) \cup ... \cup (V_{2k-1}\setminus 
V_{2k-2}).$$
Ceci termine la preuve.$\hfill\square$

\begin{cor} Soient $X$ et $Y$ des espaces polonais, $A$ un 
bor\'elien \`a coupes verticales d\'enombrables de $X \times Y$, et $k$ un entier 
naturel non nul. Alors on a les \'equivalences suivantes :\smallskip

\noindent (a) $A$ est non-$\mbox{pot}(\borone)\Leftrightarrow A$ est non-$\mbox{pot}(\boraone)
\Leftrightarrow$ la projection de $A$ sur $Y$ est non d\'enombrable.\smallskip

\noindent (b) $A$ est non-$\mbox{pot}(\check D_{2k-1}(\boraone))\Leftrightarrow A$ 
est non-$\mbox{pot}(\check D_{2k}(\boraone))\Leftrightarrow A$ est 
non-$\mbox{pot}(D_{2k-1}(\boraone)^+)$.\smallskip

\noindent (c) $A$ est non-$\mbox{pot}(D_{2k+1}(\boraone))\Leftrightarrow A$ est 
non-$\mbox{pot}(D_{2k}(\boraone))\Leftrightarrow A$ est 
non-$\mbox{pot}(D_{2k}(\boraone)^+)$.\end{cor}
  
  Ce r\'esultat, coupl\'e avec le th\'eor\`eme 2.5, donne des caract\'erisations des 
  ensembles potentiellement diff\'erences finies d'ouverts parmi les bor\'eliens \`a 
  coupes verticales d\'enombrables. Les classes $\{ \emptyset \}$, $\borone$, 
  $D_{2n+2}(\boraone)$, et $\check D_{2n+1}(\boraone)$ sont les seules de la 
  forme $\Gamma_A$, avec $A$ bor\'elien \`a coupes verticales d\'enombrables 
potentiellement diff\'erences finies d'ouverts, comme on le voit avec les 
r\'esultats pr\'ec\'edents et les ensembles suivants :
$$\begin{array}{ll}
A_{2n+2}\!\!\! & := \{(\alpha,\beta) \in 2^\omega \times 2^\omega~/~\mbox{le nombre d'entiers
   o\`u }\alpha\mbox{ et }\beta\mbox{ diff\`erent est impair }\leq 2n+1 \}\mbox{,}\cr & \cr
A_{2n+1}\!\!\! & := \{(\alpha,\beta) \in 2^\omega \times 2^\omega~/~\mbox{le nombre d'entiers
   o\`u }\alpha\mbox{ et }\beta\mbox{ diff\`erent est pair }\leq 2n \}.
\end{array}$$
V\'erifions-le pour $A_{2n+2}$, par exemple. Il est clairement $D_{2n+2}(\boraone)$, 
puisque si on pose 
$$U_p := \{(\alpha,\beta) \in 2^\omega \times 2^\omega~/~\exists~s\in 
\omega^{2n+2-p}~~i\not= j \Rightarrow s(i)\not= s(j)\mbox{ et }\forall~i<\vert s\vert ~~
\alpha(s(i)) \not= \beta(s(i)) \}\mbox{,}$$
$A_{2n+2} = \bigcup_{p<2n+2,~p~\mbox{impair}} U_p \setminus U_{p-1}$ ; 
pour voir qu'il est non-$\mbox{pot}(\check D_{2n+2}(\boraone))$, il suffit par le 
corollaire 2.9 de voir qu'il est non-$\mbox{pot}(\check D_{2n+1}(\boraone))$. On va 
pour ce faire appliquer le th\'eor\`eme 2.5. On pose 
  $f_{\emptyset}(\alpha) = \alpha$, et si $\vert t\vert \leq 2n$, 
  $f_{t^\frown m}(\alpha)(p)=f_t(\alpha)(p) \Leftrightarrow p \not= m$. Posons 
  \'egalement $\phi_m(\alpha)(p)=\alpha(p) \Leftrightarrow p \not= m$. Alors 
  $\phi_m$ est un hom\'eomorphisme, et comme $f_{t^\frown m} = \phi_m\circ f_t$, 
  on a par r\'ecurrence que les $f_t$ sont des hom\'eomorphismes de $2^\omega$ sur 
  lui-m\^eme. Par r\'ecurrence sur $\vert t\vert  \leq 2n+1$, on voit que le nombre de $p$
   tels que $\alpha(p) \not= f_t(\alpha)(p)$ a m\^eme parit\'e que $\vert t\vert $ et est 
   $\leq \vert t\vert $. On en d\'eduit que $\bigcup_{s\in \omega^{\leq 2n+1}~/~\vert s\vert ~
\mbox{paire}} C_s \subseteq A_{2n+2}$. Enfin, si $\alpha \not= \beta$ en $n_0<n_1<...<n_{2p}$, 
avec $p\leq n$, et si on pose $\alpha_0 =\alpha$, 
$\alpha_{l+1} = \phi_{n_l}(\alpha_l)$ si $l\leq 2p$, $\alpha_{2p+1} = \beta$, et 
$s = (n_0,...,n_{2p})$, on v\'erifie que si $q\leq 2p+1$, 
$\alpha_q = f_{s \lceil q}(\alpha)$ ; d'o\`u $(\alpha,\beta) \in C_{s\lceil 2p}$.
 
\vfill\eject 
 
 Dans [SR], il est d\'emontr\'e le r\'esultat suivant, d\^u \`a Hurewicz :

\begin{thm} Soit $X$ un espace polonais, et $A$ un bor\'elien de $X$. 
Alors $A$ est non-$\bormtwo$ si et seulement s'il existe $E$ d\'enombrable 
sans point isol\'e tel que $\overline{E}\setminus E \approx \omega^\omega$ et $E = A \cap 
\overline{E}$.\end{thm}

\begin{thm} Soient $X$ et $Y$ des espaces polonais, 
$A$ un bor\'elien de $X\times Y$ \`a coupes verticales 
d\'enombrables. Alors $A$ est non-$\mbox{pot}(\bormtwo)$ si et 
seulement s'il existe des espaces polonais $Z$ et $T$ parfaits de dimension
 0 non vides, des injections continues $u$ et 
$v$, une suite $(A_n)$ d'ouverts denses de $Z$, une suite $(f_n)$ 
d'applications continues et ouvertes de $A_n$ dans $T$, tels que pour tout 
$x$ dans $\bigcap_{n\in\omega} A_n$, l'ensemble $E_x := \{ 
f_n(x)~/~n\in\omega \}$ soit sans point isol\'e, 
$\overline{E_x} \setminus E_x \approx \omega^\omega$, et aussi 
$E_x = (u\times v)^{-1}(A)_x \cap \overline{E_x}$.\end{thm}

\noindent\bf D\'emonstration.\rm\ Supposons que $A$ soit $\mbox{pot}(\bormtwo)$, en raisonnant par 
l'absurde. Alors il existe un $G_{\delta}$ dense $K$ de $T$ tel que les coupes 
verticales de $(Z\times K) \cap (u\times v)^{-1}(A)$ soient $G_{\delta}$ dans 
$K$, donc dans $T$. Par ailleurs, 
$\bigcap_{n\in\omega} A_n \cap \bigcap_{n\in\omega} f_n^{-1}(K)$ est $G_\delta$ 
dense de $Z$, puisque les fonctions $f_n$ sont continues et ouvertes. Donc on 
trouve $x$ dans ce $G_\delta$, et $E_x$ est polonais, puisque c'est une 
coupe verticale de $(Z\times K) \cap (u\times v)^{-1}(A)$ intersect\'ee avec 
$\overline{E_x}$. De plus $E_x$ est d\'enombrable, non vide et sans point 
isol\'e, ce qui est contradictoire.\bigskip

 Inversement, on peut supposer, pour simplifier l'\'ecriture, que $X$ 
et $Y$ sont r\'ecursivement pr\'esen-t\'es, et que $A$ est 
$\Borel$-r\'eunion de graphes $\Borel$. D\'esignons par $W^X$ un 
ensemble $\Ca \subseteq \omega$ de codes pour les $\Borel$ de $X$, et 
par $C^X \subseteq \omega\times X$ un ensemble $\Ca$ dont les sections 
aux points de $W^X$ d\'ecrivent les $\Borel$ de $X$, et tel que la relation 
($n\in W^X$ et $(n,x)\notin C^X$) soit $\Ca$ (cf [Lo1]). Soit \'egalement $W\subseteq W^{X\times Y}$ un 
ensemble $\Ca$ de codes pour les $\Borel \cap \mbox{pot}(\boratwo)$ de $X\times Y$ (dont l'existence 
est d\'emontr\'ee dans [Lo2]). Posons 
$$H := \cup\ \{\ (E\times F)\setminus A~/~E,F\in \Borel\mbox{ ~et~ }(E\times F)\setminus A\in \mbox{pot}(\boratwo)\ \}.$$ 
On a 
$$\begin{array}{ll}
H(x,y)~\ \Leftrightarrow
& \exists~n\in W~\exists~m~:~
(m)_0 \in W^X ~\mbox{et}~(m)_1 \in W^Y~\mbox{et}~\forall~z~~\forall~t \cr 
& [(n\in W^{X\times Y}~\mbox{et}~(n,z,t)\notin C^{X\times Y})~\mbox{ou}\cr
& \{((m)_0,z)\in C^X~\mbox{et}~((m)_1,t)\in C^Y~\mbox{et }(z,t)\notin A\}]~\mbox{et}\cr 
& [((m)_0\in W^X~\mbox{et}~((m)_0,z)\notin C^X)\mbox{ ou}~
((m)_1\in W^Y~\mbox{et}~((m)_1,t)\notin C^Y)~\mbox{ou}\cr 
& (z,t)\in A~\mbox{ou}~(n,z,t)\in C^{X\times Y}]\cr & \mbox{et}~
((m)_0,x)\in C^X~\mbox{et}~((m)_1,y)\in C^Y~\mbox{et}~(x,y)\notin A.
\end{array}$$
Donc $H$ est $\Ca$. Posons $N := \check A \cap \check H$ ; $N$ est $\Ana$ et 
$G_{\delta}$ de $X\times Y$ muni de la topologie 
${\it\Delta}_X \times {\it\Delta}_Y$ (cf [Lo2]). Posons $D_X := \{ x\in X~/~x\notin \Borel \}$, 
$\Omega_X := \{ x\in X~/~\omega_1^x = \omega_1^{\mbox{CK}} \}$, et 
$Z_0 := \Omega_X \cap D_X$, $T := \Omega_Y \cap D_Y$. On munit $Z_0$ 
(resp. $T$) de la restriction de la topologie de Gandy-Harrington sur $X$ 
(resp. $Y$), de sorte que $Z_0$ et $T$ sont polonais parfaits de dimension 0. 
 On montre maintenant la propri\'et\'e qui sera la clef de la construction \`a venir :\bigskip
 
 Si $U$ (resp. $V$) est $\Ana$ et inclus dans $Z_0$ (resp. $T$), et si 
$N\cap (U\times V) \not= \emptyset$, alors l'ensemble 
$A\cap \overline{N\cap (U\times V)}^{Z_0\times T}$ est un sous-$\Ana$ 
non vide de $U\times V$.

\vfill\eject

 En effet, $U$ (resp. $V$) est ouvert-ferm\'e de $Z_0$ (resp. $T$), donc 
$\overline{N\cap (U\times V)}^{Z_0\times T} \subseteq U\times V$ ; $A$, $Z_0$, 
$T$, $N$, $U$, et $V$ sont $\Ana$, donc $A\cap \overline{N\cap (U\times V)}^{Z_0\times T}$ 
aussi, l'adh\'erence d'un $\Ana$ pour le produit des topologies de Gandy-Harrington 
restant $\Ana$, par double application du th\'eor\`eme de s\'eparation. Posons 
$O := N\cap (U\times V)$. Supposons que $A\cap \overline{O}^{{\it\Sigma}_X\times 
{\it\Sigma}_Y} = \emptyset$. Alors, par double application du th\'eor\`eme de s\'eparation, 
$A\cap \overline{O}^{{\it\Delta}_X\times {\it\Delta}_Y} = \emptyset$, donc on a la triple inclusion 
$$(U\times V)\setminus A\subseteq O \cup H \subseteq
\overline{O}^{{\it\Delta}_X\times {\it\Delta}_Y}\cup H \subseteq \check A.$$
Donc $(U\times V)\setminus A$ et $A$ sont deux $\Ana$ s\'eparables par un 
$\mbox{pot}(\boratwo)$ ; ils peuvent par cons\'equent \^etre s\'epar\'es par un ensemble 
$K\in \Borel \cap \mbox{pot}(\boratwo)$ (cf [Lo2]). On a $U\times V \subseteq 
K\cup A$, donc on peut trouver $\cal U$ et $\cal V$, deux $\Borel$ tels que 
$U\times V\subseteq \cal U\times \cal V$ $ \subseteq K\cup A$. On a donc 
$(\cal U\times \cal V)$$\setminus A\subseteq H$, puis $O \subseteq H\setminus H = \emptyset$, 
ce qui est absurde. On a donc que $A\cap \overline{O}^{{\it\Sigma}_X\times 
{\it\Sigma}_Y} \not= \emptyset$. Or $O\subseteq D_X\times D_Y$, donc $\overline{O}^{{\it\Sigma}_X\times 
{\it\Sigma}_Y} \subseteq D_X\times D_Y$. Donc $(D_X\times D_Y)\cap A\cap \overline{O}^{{\it\Sigma}_X\times {\it\Sigma}_Y} \not= \emptyset$ et est $\Ana$, donc rencontre $\Omega_{X\times Y} \subseteq \Omega_X \times \Omega_Y$, donc $(Z_0\times T)\cap A\cap \overline{O}^{{\it\Sigma}_X\times 
{\it\Sigma}_Y} \not= \emptyset$.\bigskip

 Soit $d_1$ (resp. $d_2$) une distance $\leq 1/2$ qui rende $Z_0$ (resp. $T$) 
complet, et $d$ la distance sur $Z\times T$ d\'efinie par 
$d((x,y),(z,t)) := d_1(x,z) + d_2(y,t)$. Soit $d'$ une distance $\leq 1$ sur $N\cap (Z_0\times T)$, qui 
le rende complet (c'est un $G_\delta$ de $Z_0\times T$, donc un espace polonais).
On v\'erifie sans peine que si on pose $d'_x(y,t) := d'((x,y),(x,t))$, alors 
$(N_x\cap T,d'_x)$ est m\'etrique complet. On pose, si $s$ est une suite finie 
d'entiers non nuls, $\nu (s) := \Sigma_{i<\vert s\vert }~s(i)$ ~($\nu (\emptyset) = 0$).\bigskip 

On construit, par r\'ecurrence sur $\vert s\vert $, des ouverts non vides de $Z_0$, 
${\cal A}_s$, 
des $G_\delta$ denses ${\cal G}_s$ de ${\cal A}_s$, des applications continues et ouvertes 
$g_s:{\cal A}_s \rightarrow T$, des ouverts \`a coupes verticales ouvertes-ferm\'ees 
de $Z_0 \times T$, $\omega_s$, et des ouverts de $N\cap (Z_0\times T)$, $G_s$, tels que 
$$\begin{array}{ll} 
& (a)~\mbox{Gr}(g_s)\subseteq A\cap\omega_s \cap \overline{G_s}^{Z_0\times T} \cr 
& (b)~\omega_{s^\frown n} \subseteq \omega_s,~\omega_{s^\frown n}\cap\omega_{s^\frown m}
= \emptyset~\mbox{si}~n\not= m \cr 
& (c)~\forall~x\in Z_0~~\overline{G_{s^\frown n}^x}^T\cap N_x \subseteq G_s^x 
\subseteq \omega_s^x \cr 
& (d)~{\cal A}_s\Delta {\cal A}_{s^\frown n}~\mbox{est~maigre~et}~\forall~x\in 
{\cal A}_{s^\frown n}\cap{\cal G}_s~~d_2(g_s(x),g_{s^\frown n}(x)) \leq 2^{1-\nu (s^\frown n)} \cr
& (e)~\forall~x\in Z_0~~\delta_2(\omega_s^x) \leq 2^{-\nu (s)} \cr 
& (f)~\forall~x\in Z_0~~\delta'_x(G_s^x) \leq 2^{-\vert s\vert }
\end{array}$$
On pose $\omega_\emptyset := Z_0\times T$, $G_\emptyset := (Z_0\times T)\cap N$. 
L'ensemble $G_{\emptyset}$ est non vide. En effet, $N$ 
rencontre $D_X\times D_Y$, sinon $N$ serait $\mbox{pot}(\borone)$ et 
$A\cup H$ aussi, donc $A = (A\cup H)\cap \check H$ serait 
$\mbox{pot}(\bormtwo)$, ce qui est exclus. Donc $N\cap (D_X\times D_Y)$, 
qui est $\Ana$, rencontre $\Omega_{X\times Y} \subseteq \Omega_X 
\times \Omega_Y$ et $G_\emptyset \not= \emptyset$. Par la propri\'et\'e-clef, 
$A\cap \overline{G_\emptyset}^{Z_0\times T}$ est un sous-$\Ana$ non vide de 
$\omega_\emptyset$, donc $A\cap \omega_\emptyset \cap \overline{G_\emptyset}^{Z_0\times T}$ 
est r\'eunion de graphes $\Ana$ dont l'un, disons $\mbox{Gr}(g_\emptyset)$, est non vide. 
Il reste \`a appeler ${\cal A}_\emptyset$ le domaine de $g_\emptyset$ pour terminer la 
construction au premier cran. Admettons donc avoir construit ${\cal A}_s$, $g_s$, 
$\omega_s$, $G_s$ pour $\vert s\vert  \leq p$ et ${\cal G}_s$ pour $\vert s\vert  < p$ v\'erifiant (a)-(f), 
et soit $s$ de longueur $p$.

\vfill\eject

 On commence par construire, par r\'ecurrence sur $n$, les $\omega_{s^\frown n}$, 
une suite d\'ecroissante $(E_n)$ de $G_\delta$ denses de ${\cal A}_s$, des applications 
continues $h_n : E_n \rightarrow T$, et des ouverts $V_n$ de $Z_0\times T$ tels que 
$$\begin{array}{ll} 
 & (1)~\mbox{Gr}(h_n) \subseteq \omega_{s^\frown n}\cap G_s \cr
 & (2)~\omega_{s^\frown n} \subseteq \omega_s~\mbox{et}~
  \omega_{s^\frown n}\cap\omega_{s^\frown m} = \emptyset~\mbox{si}~n\not= m \cr 
 & (3)~\forall~x\in E_n~~d_2(h_n(x),g_s(x)) < 2^{-\nu (s^\frown n)} \cr 
 & (4)~\forall~x\in Z_0~~\delta_2(\omega_{s^\frown n}^x) < 2^{-\nu (s^\frown n)} \cr 
 & (5)~\mbox{Gr}(g_s\lceil E_n) \subseteq V_n \subseteq~^c(\bigcup_{1\leq m\leq n} 
\omega_{s^\frown m})
\end{array}$$ 
Admettons avoir construit ces objets pour $1\leq m\leq n$, ce qui est fait pour $n=0$. 
Soit $x$ un \'el\'ement de ${\cal A}_s$. Comme $g_s$ est continue en $x$, 
on trouve un voisinage ouvert $U_x$ de $x$ contenu dans ${\cal A}_s$ tel que 
si $z\in U_x$, $d_2(g_s(x),g_s(z)) < \varepsilon := 2^{-\nu (s^\frown (n+1))-1}$. De 
telle fa\c con que si $(z,t)\in U_x\times {\cal B}(g_s(x),\varepsilon)$, on a $d_2(t,g_s(z)) < 
2^{-\nu (s^\frown (n+1))}$. Posons ${\cal U}_{n+1} := V_n\cap \bigcup_{x\in {\cal A}_s} 
U_x\times {\cal B}(g_s(x),\varepsilon)$, et soit $O_n$ un ouvert dense de ${\cal A}_s$ 
contenant $E_n$ tel que $\mbox{Gr}(g_s)\cap V_n = \mbox{Gr}(g_s\lceil O_n)$. Soit $k_{n+1}$ une fonction 
Baire-mesurable uniformisant ${\cal U}_{n+1}\cap \omega_s \cap G_s$ sur sa projection $\Pi$ 
qui est un ouvert de $Z_0$, et soit $F_{n+1}$ un $G_\delta$ dense de $\Pi$ sur 
lequel la restriction de $k_{n+1}$ est continue. Alors $O_n \setminus F_{n+1}$ 
est maigre car $O_n\cap \Pi$ est ouvert dense de $O_n$ ; en effet, si $U$ est ouvert 
non vide de $O_n$ et $x \in U$, $(x,g_s(x))\in (U\times T)\cap \omega_s \cap {\cal U}_{n+1}$. 
Donc, comme $\mbox{Gr}(g_s) \subseteq \overline{G_s}^{Z_0\times T}$, on trouve $(z,t)$ dans 
$G_s\cap (U\times T)\cap \omega_s \cap {\cal U}_{n+1}$, et $z\in U\cap \Pi$. De 
tout ceci r\'esulte que ${\cal A}_s\setminus F_{n+1}$ est maigre, donc $E_{n+1} := 
F_{n+1}\cap E_n$ est $G_\delta$ dense de ${\cal A}_s$ et $h_{n+1} := k_{n+1}\lceil E_{n+1}$ 
est continue. Si $x\in E_{n+1}$, $(x,h_{n+1}(x))\in {\cal U}_{n+1}$, donc 
$$d_2(h_{n+1}(x),g_s(x)) < 2^{-\nu (s^\frown (n+1))}.$$ 
Dans $E_{n+1}\times T$, $\mbox{Gr}(h_{n+1})$ et $\mbox{Gr}(g_s\lceil E_{n+1})$ sont deux ferm\'es disjoints, donc s\'eparables par un 
ouvert-ferm\'e $\theta$. Par la propri\'et\'e de r\'eduction des ouverts, on trouve donc 2 ouverts disjoints de $Z_0 \times T$, $\cal T$ et $\cal W$ tels que 
$\theta = {\cal T}\cap (E_{n+1}\times T)$, et $(E_{n+1}\times T)\setminus \theta = {\cal W}\cap (E_{n+1}\times T)$. 
Comme $\mbox{Gr}(h_{n+1}) \subseteq \cal T$, on trouve pour chaque $x$ dans $E_{n+1}$ un rectangle ouvert-ferm\'e 
$V_x\times W_x$ tel que 
$$(x,h_{n+1}(x))\in V_x\times W_x \subseteq {\cal T}\capÊ\omega_s \cap 
V_n\mbox{, }$$
$\delta_2(W_x) < 2^{-\nu (s^\frown (n+1))}$, et $E_{n+1}\cap V_x \subseteq h_{n+1}^{-1}(W_x)$.  
On a $\bigcup_{x\in E_{n+1}} V_x\times W_x = \bigcup_{p\in\omega} V_p\times W_p$, 
par Lindel\"of. R\'eduisons la suite $(V_p)$ en $(V'_p)$ ; on peut alors poser 
$\omega_{s^\frown (n+1)} := \bigcup_{p\in\omega} V'_p\times W_p$ et $V_{n+1} 
:= V_n \cap \cal W$, comme on le v\'erifie facilement.\bigskip

 Revenons \`a la construction principale ; on a d\'ej\`a assur\'e (b) et (e). Si 
$x\in E_n$, on peut trouver un ouvert-ferm\'e de $N$, ${\cal U}_x := ({\cal V}_x 
\times {\cal W}_x)\cap N$, de diam\`etre au plus $2^{-\vert s\vert -1}$ pour $d'$, v\'erifiant 
$${\cal V}_x\cap E_n \subseteq h_n^{-1}({\cal W}_x)\mbox{,}$$ 
et si $Z_s$ est ouvert de $Z_0\times T$ tel que $G_s = Z_s \cap N$, on ait 
$$(x,h_n(x)) \in {\cal V}_x\times {\cal W}_x \subseteq \omega_{s^\frown n} 
\cap Z_s \cap ({\cal A}_s\times T).$$ 
On a donc que $\mbox{Gr}(h_n) \subseteq \bigcup_{x\in E_n} {\cal U}_x = 
\bigcup_{m\in \omega} {\cal U}_m$, par Lindel\"of. On r\'eduit la suite $({\cal V}_m)$ 
en $({\cal V}'_m)$, et on pose $G_{s^\frown n} := \bigcup_{m\in\omega} ({\cal V}'_m\times 
{\cal W}_m)\cap N$, et les conditions (c) et (f) sont r\'ealis\'ees, puisque $G_{s^\frown n}^x$ est vide ou 
l'un des ${\cal W}_m\cap N_x$, donc est ouvert-ferm\'e de $N_x$. De plus, on a 
$G_{s^\frown n} \supseteq \mbox{Gr}(h_n)$.

\vfill\eject

 A l'aide d'une nouvelle num\'erotation, on 
peut \'ecrire $G_{s^\frown n} = \bigcup_{p\in\omega} 
(L_p\times M_p)\cap N$, $L_p$ et $M_p$ \'etant deux $\Ana$. Posons 
$E_{n,p} := A\cap \overline{N\cap (L_p\times M_p)}^{Z_0\times T}\cap 
\omega_{s^\frown n}$ ; $E_{n,p}$ est r\'eunion 
de graphes $\Ana$. On peut donc trouver 
une suite de graphes $\Ana$, contenus dans $E_{n,p}$, dont les domaines sont deux \`a deux 
disjoints, et tels que la r\'eunion de ces domaines, disons $O_{n,p}$ forme un ouvert 
tel que $\bigcup_{q\leq p} O_{n,q}$ soit dense dans 
$\bigcup_{q\leq p} \Pi_X''E_{n,p}$ (une telle construction est possible avec le lemme 
de Zorn, par exemple). On fait ceci par r\'ecurrence sur $p$, de sorte que 
${\cal A}_{s^\frown n} :=  \bigcup_{p\in\omega} O_{n,p}$ est un ouvert dense de 
$\bigcup_{p\in\omega} \Pi_X''E_{n,p}$. On appelle $g_{s^\frown n}$ le recollement 
des fonctions d\'efinies sur les $O_{n,p}$, et la condition (a) est satisfaite, 
ainsi que la seconde partie de la condition (d), si on pose ${\cal G}_s := 
\bigcap_{n\geq 1} E_n$, par (4). Reste \`a voir que ${\cal A}_s\Delta {\cal A}_{s^\frown n}$ 
est maigre pour clore la construction. Quitte \`a remplacer ${\cal A}_{s^\frown n}$ par 
${\cal A}_s \cap {\cal A}_{s^\frown n}$, il suffit de voir que 
${\cal A}_s\setminus {\cal A}_{s^\frown n}$ est maigre. Posons 
$\Pi_p := \Pi_X''(L_p\times M_p)\cap N$ et 
$\Pi'_p := \Pi_X''E_{n,p}$. On a 
$$\begin{array}{ll}
{\cal A}_s\setminus {\cal A}_{s^\frown n}\!\!\! 
& \subseteq {\cal A}_s
\setminus (\bigcup_{p\in\omega} \Pi'_p) \cup \bigcup_{p\in\omega} \Pi'_p
\setminus {\cal A}_{s^\frown n} \cr & \subseteq {\cal A}_s\setminus E_n 
\cup E_n\setminus (\bigcup_{p\in\omega} \Pi'_p) \cup \bigcup_{p\in\omega} 
\Pi'_p\setminus {\cal A}_{s^\frown n} \cr 
& \subseteq {\cal A}_s\setminus E_n \cup 
\bigcup_{p\in\omega} \Pi_p\setminus \Pi'_p \cup \bigcup_{p\in\omega} 
\Pi'_p\setminus {\cal A}_{s^\frown n}.
\end{array}$$ 
Il suffit donc de voir que $\Pi_p\cap\Pi'_p$ est ouvert dense de $\Pi_p$ ; mais ceci est une cons\'equence facile de la propri\'et\'e clef. \bigskip

 On d\'efinit maintenant les objets recherch\'es : $Z:={\cal A}_\emptyset$, $u$ 
et $v$ sont les injections canoniques. Soit $\varphi : \omega \rightarrow (\omega\setminus\{
0\})^{<\omega}\times \omega$ bijective, et $(U_m^s)_m$ une suite d'ouverts denses de ${\cal A}_s$ 
telle que ${\cal G}_s = \bigcap_{m\in\omega} U_m^s$. On pose $A_n := {\cal A}_\emptyset 
\cap U_{\varphi_1(n)}^{\varphi_0(n)}$ et $f_n := g_{\varphi_0(n)}\lceil A_n$. L'ouvert $A_n$ 
est dense dans $Z$ par la premi\`ere partie de (d), et $E_x$ est sans point isol\'e 
par la seconde partie de (d) (on a l'\'egalit\'e $\bigcap_{n\in\omega} A_n = \bigcap_{s\in
(\omega\setminus\{ 0\})^{<\omega}} {\cal G}_s$).\bigskip 

 La fin de la preuve est semblable \`a la preuve du th\'eor\`eme d'Hurewicz dans [SR]. 
Posons, si $x\in \bigcap_{n\in\omega} A_n$, $M_k := \{ g_s(x)~/~\vert s\vert  \leq k\}$. 
Alors $M_k$ est ferm\'e dans $T$, par r\'ecurrence sur k : si $M_k^\varepsilon := \{ y\in 
T~/~d_2(y,M_k)\leq \varepsilon\}$, on a $M_{k+1} = \bigcap_{\varepsilon > 0} [M_k^\varepsilon \cup 
(M_{k+1}\setminus M_k^\varepsilon)]$, et $M_{k+1}\setminus M_k^\varepsilon$ est fini, par (d).\bigskip

 Si $k\in\omega$ et $y\in \overline{E_x}\setminus E_x$, $y\notin M_k$, donc il 
existe $\varepsilon > 0$ tel que $y\notin M_k^\varepsilon$. Si $s\in (\omega\setminus\{0\})^k$ et $(x,t)\in 
\omega_s$, $d_2(t,M_k)\leq d_2(t,g_s(x))\leq \delta_2(\omega_s^x)\leq 2^{-\nu (s)} 
\leq \varepsilon$ d\`es que $\nu (s)\geq k_0$, donc $\omega_s^x \subseteq M_k^\varepsilon$, 
sauf pour un nombre fini de $s$ dans $(\omega\setminus\{0\})^k$. Donc si 
${\cal H} := \{s\in (\omega\setminus\{0\})^k~/~\omega_s^x \not\subseteq M_k^\varepsilon\}$, 
on a $E_x = M_k \cup \{ g_s(x)~/~\vert s\vert >k\} \subseteq M_k \cup \bigcup_{\vert t\vert >k} 
\omega_t^x \subseteq M_k \cup \bigcup_{\vert s\vert  = k} \omega_s^x \subseteq M_k^\varepsilon 
\cup \bigcup_{s\in {\cal H}} \omega_s^x$ et $\overline{E_x} \subseteq M_k^\varepsilon 
\cup \bigcup_{s\in {\cal H}} \omega_s^x \subseteq M_k^\varepsilon 
\cup \bigcup_{\vert s\vert  = k} \omega_s^x$. Donc on trouve une unique 
suite $\sigma$ dans $(\omega\setminus\{0\})^\omega$ telle que $y\in \bigcap_{s
\prec \sigma} \omega_s^x$. La suite d\'ecroissante de ferm\'es non vides 
dont les diam\`etres tendent vers 0 de $(N_x\cap T, d'_x)$, $(\overline{G_s^x}
\cap N_x)_{s\prec \sigma}$, converge vers $\xi\in N_x\cap T$, et $\{\xi\} = 
\bigcap_{s\prec\sigma} G_s^x$. 
D'o\`u $\xi \in \bigcap_{s\prec\sigma} \omega_s^x = \{y\}$ et $(x,y)\in 
N\subseteq \check A$.~Posons 
$$h:\left\{\!\!
\begin{array}{ll}
(\omega\setminus\{0\})^\omega\!\! 
& \rightarrow T \cr \sigma\!\! 
& \mapsto \bigcap_{s\prec\sigma} 
\overline{G_s^x}
\end{array}\right.$$
La fonction $h$ est bien d\'efinie, \`a valeurs dans 
${N_x\subseteq \check E_x}$. Si $\sigma \in (\omega\setminus\{0\})^\omega$, 
$h(\sigma) \in \bigcap_{s\prec\sigma} \omega_s^x$ et 
$\omega_s^x\cap E_x \not= \emptyset$ ; par cons\'equent, $h$ est \`a valeurs dans 
$\overline{E_x}$. On a vu que $\overline{E_x}\setminus E_x \subseteq 
h''(\omega\setminus\{0\})^\omega$. Par (b), $h$ est injective, et par (e), $h$ est 
continue. Comme $h''N_s = \omega_s^x \cap (\overline{E_x}\setminus E_x)$, $h$ 
est bien un hom\'eomorphisme de $(\omega\setminus\{0\})^\omega$ sur 
$\overline{E_x}\setminus E_x$.$\hfill\square$

\vfill\eject

 Avec ce r\'esultat et le th\'eor\`eme 2.3, on a, pour chaque 
classe de Wadge $\Gamma$ non auto-duale, un test pour dire 
si un bor\'elien \`a coupes d\'enombrables est $\mbox{pot}(\Gamma)$ (en 
effet, un tel bor\'elien est $\mbox{pot}(F_\sigma)$).

\section{$\!\!\!\!\!\!$ R\'ef\'erences.}

\noindent [D-SR]\ \ G. Debs et J. Saint Raymond,~\it S\'elections
bor\'eliennes injectives,~\rm Amer. J. Math.~111 (1989),
519-534

\noindent [GM]\ \ S. Graf and R. D. Mauldin,~\it Measurable one-to-one
selections and transition kernels,~\rm Amer. J. Math.~107 (1985), 407-425

\noindent [HKL]\ \ L. A. Harrington, A. S. Kechris et A. Louveau,~\it A Glimm-Effros 
dichotomy for Borel equivalence relations,~\rm J. Amer. Math. Soc.~3 (1990), 903-928

\noindent [Ke]\ \ A. S. Kechris,~\it Measure and category in effective 
descriptive set theory,~\rm Ann. of Math. Logic,~5 (1973), 337-384

\noindent [Ku]\ \  K. Kuratowski,~\it Topology,~\rm Vol. 1, Academic Press, 1966

\noindent [Le1]\ \ D. Lecomte,~\it Classes de Wadge potentielles et th\'eor\`emes d'uniformisation 
partielle,~\rm Fund. Math.~143 (1993), 231-258

\noindent [Le2]\ \ D. Lecomte,~\it Classes de Wadge potentielles des bor\'eliens \`a coupes d\'enombrables,\ \rm C. R. Math. Acad. Sci. Paris\ 317, S\'erie 1 (1993), 1045-1048

\noindent [Lo1]\ \ A. Louveau,~\it Livre \`a para\^\i tre\rm

\noindent [Lo2]\ \ A. Louveau,~\it Ensembles analytiques et bor\'eliens dans les espaces produit,~\rm Ast\'erisque (S. M. F.) 78 (1980)

\noindent [Lo-SR]\ \ A. Louveau and J. Saint Raymond,~\it Borel classes and closed games : Wadge-type and Hurewicz-type results,~\rm Trans. A. M. S.~304 (1987), 431-467

\noindent [Ma]\ \ R. D. Mauldin,~\it One-to-one selections, marriage theorems,~\rm Amer. J. Math.~104 (1982), 823-828

\noindent [Mo]\ \ Y. N. Moschovakis,~\it Descriptive set theory,~\rm North-Holland, 1980

\noindent [O]\ \ J. C. Oxtoby,~\it Measure and category,~\rm Springer-Verlag, 1971

\noindent [S]\ \ J. R. Steel,~\it Analytic sets and Borel isomorphisms,~\rm 
Fund. Math.~108 (1980), 83-88

\noindent [SR]\ \ J. Saint Raymond,~\it La structure bor\'elienne d'Effros est-elle standard ?,~\rm 
Fund. Math.~100 (1978), 201-210

\noindent [W]\ \ W. W. Wadge,~\it Thesis,~\rm Berkeley (1984) 

\end{document}